\documentclass[12pt,leqno]{article}
\usepackage{amssymb}
\usepackage[mathscr]{eucal}
\usepackage{amsmath,amssymb,latexsym,theorem,bbm}
\usepackage{color,url}
\usepackage[hyperfootnotes=false]{hyperref}

\newcommand{\SC}{\scriptstyle}

\newcommand{\CC}{\mathsf{C}}
\newcommand{\DD}{\mathsf{D}}

\newcommand{\NN}{\mathbb{N}}

\newcommand{\RR}{\mathbb{R}}

\newcommand{\ZZ}{\mathbb{Z}}

\newcommand{\bA}{{\boldsymbol{A}}}

\newcommand{\bb}{{\boldsymbol{b}}}
\newcommand{\bB}{{\boldsymbol{B}}}

\newcommand{\bC}{{\boldsymbol{C}}}

\newcommand{\bI}{{\boldsymbol{I}}}

\newcommand{\bM}{{\boldsymbol{M}}}
\newcommand{\bp}{{\boldsymbol{p}}}
\newcommand{\bP}{{\boldsymbol{P}}}

\newcommand{\bR}{{\boldsymbol{R}}}

\newcommand{\bS}{{\boldsymbol{S}}}

\newcommand{\bT}{{\boldsymbol{T}}}
\newcommand{\bu}{{\boldsymbol{u}}}

\newcommand{\bV}{{\boldsymbol{V}}}

\newcommand{\bx}{{\boldsymbol{x}}}
\newcommand{\bX}{{\boldsymbol{X}}}

\newcommand{\bU}{{\boldsymbol{U}}}

\newcommand{\bgamma}{{\boldsymbol{\gamma}}}
\newcommand{\bxi}{{\boldsymbol{\xi}}}

\newcommand{\bbeta}{{\boldsymbol{\beta}}}
\newcommand{\bvare}{{\boldsymbol{\vare}}}

\newcommand{\bzero}{{\boldsymbol{0}}}

\newcommand{\cA}{{\mathcal A}}
\newcommand{\cB}{{\mathcal B}}

\newcommand{\cD}{{\mathcal D}}

\newcommand{\cF}{{\mathcal F}}
\newcommand{\cH}{{\mathcal H}}

\newcommand{\cM}{{\mathcal M}}
\newcommand{\bcM}{\boldsymbol{\cM}}

\newcommand{\cP}{{\mathcal P}}
\newcommand{\bcP}{\boldsymbol{\cP}}

\newcommand{\cS}{{\mathcal S}}
\newcommand{\cU}{{\mathcal U}}

\newcommand{\bcU}{\boldsymbol{\cU}}

\newcommand{\cX}{{\mathcal X}}
\newcommand{\bcX}{\boldsymbol{\cX}}

\newcommand{\cY}{{\mathcal Y}}
\newcommand{\cZ}{{\mathcal Z}}
\newcommand{\cW}{{\mathcal W}}
\newcommand{\bcW}{\boldsymbol{\cW}}

\newcommand{\bcY}{\boldsymbol{\cY}}

\newcommand{\bcZ}{\boldsymbol{\cZ}}

\newcommand{\dd}{\mathrm{d}}

\newcommand{\slu}{{\SC\mathrm{lu}}}

\newcommand{\EE}{\operatorname{\mathbb{E}}}
\newcommand{\PP}{\operatorname{\mathbb{P}}}
\newcommand{\OO}{\operatorname{O}}

\newcommand{\var}{\operatorname{Var}}
\newcommand{\cov}{\operatorname{Cov}}

\newcommand{\vare}{\varepsilon}

\renewcommand{\mid}{\,|\,}
\newcommand{\bmid}{\,\big|\,}

\renewcommand{\leq}{\leqslant}
\renewcommand{\geq}{\geqslant}

\newcommand{\skor}[1]{\stackrel{\cS_{#1}}{\longrightarrow}}
\newcommand{\stoch}{\stackrel{\PP}{\longrightarrow}}
\newcommand{\distr}{\stackrel{\cD}{\longrightarrow}}
\newcommand{\distrf}{\stackrel{\cD_f}{\longrightarrow}}
\newcommand{\distre}{\stackrel{\cD}{=}}

\newcommand{\lu}{\stackrel{\slu}{\longrightarrow}}

\newcommand{\bbone}{\mathbbm{1}}
\newcommand{\ns}{{\lfloor ns\rfloor}}
\newcommand{\nt}{{\lfloor nt\rfloor}}
\newcommand{\nT}{{\lfloor nT\rfloor}}
\newcommand{\nr}{{\lfloor nr\rfloor}}
\newcommand{\proofend}{\hfill\mbox{$\Box$}}

\numberwithin{equation}{section}

\theoremstyle{change} \theorembodyfont{\em}
\newtheorem{Lem}{Lemma.}[section]
\newtheorem{Thm}[Lem]{Theorem.}

\newtheorem{Cor}[Lem]{Corollary.}
\newtheorem{Def}[Lem]{Definition.}

\theorembodyfont{\rm}
\newtheorem{Rem}[Lem]{Remark.}

\begin{document}

\begin{center}
{\bfseries\Large Asymptotic behavior of some strongly critical\\[-1.5mm]
 decomposable multi-type Galton--Watson processes\\[1.5mm] with immigration}

\vspace*{3mm}

 {\sc\large
  M\'aty\'as $\text{Barczy}^{*,\diamond}$ \text{and}
  \ D\'aniel $\text{Bezd\'any}^{**}$ }

\end{center}

\vskip0.2cm

\noindent
 * HUN-REN--SZTE Analysis and Applications Research Group,
   Bolyai Institute, University of Szeged,
   Aradi v\'ertan\'uk tere 1, H--6720 Szeged, Hungary.

\noindent
 ** Bolyai Institute, University of Szeged,
    Aradi v\'ertan\'uk tere 1, H--6720 Szeged, Hungary.

\noindent E-mails: barczy@math.u-szeged.hu (M. Barczy),
                   bezdany@server.math.u-szeged.hu (D. Bezd\'any).

\noindent $\diamond$ Corresponding author.

\renewcommand{\thefootnote}{}
\footnote{\textit{2020 Mathematics Subject Classifications\/}: 60J80, 60F17. }
\footnote{\textit{Key words and phrases\/}:
  multi-type Galton-Watson process with immigration, strongly critical, decomposable, asymptotic behavior, squared Bessel process, iterated integral process.}
\footnote{D\'aniel Bezd\'any was supported by the 2025-2026 grant of the M\'oricz Doktorandusz Alap\'itv\'any.}

\vspace*{0.05cm}

\begin{abstract}
We study the asymptotic behavior of a critical decomposable multi-type Galton--Watson process with immigration
 whose offspring mean matrix  is triangular with diagonal entries 1.
Under second-order moment conditions on the offspring and immigration distributions,
 we establish a functional limit theorem for the sequence of appropriately scaled processes,
 and describe the limit process as independent squared Bessel processes and 
 some linear combinations of their iterated integral processes.
Our work completes and extends some results of Foster and Ney (1978) for
 some strongly critical decomposable multi-type Galton--Watson processes with immigration.
\end{abstract}

\section{Introduction}\label{Intro}

Deriving conditional and unconditional limit theorems for branching processes has a long tradition and it is an active area of research.
In this paper, we investigate the asymptotic behavior of some critical decomposable multi-type Galton-Watson processes with immigration (GWI processes).
For an overview of the history of results on the asymptotic behavior of single- and multi-type critical GWI processes,
 see the introduction of Barczy et al.\ \cite{BarBezPap2}.
In what follows, we only mention and summarize those earlier results that are directly connected to our present paper.

Foster and Ney \cite[Theorems 4 and 5]{FosNey2} proved limit theorems for some special strongly critical decomposable multi-type GWI processes.
They assumed that the offspring mean matrix of a multi-type GWI process is (upper) block triangular such that
 the submatrices in its diagonal are primitive with spectral radius 1, and that the (upper)
first superdiagonal submatrices have at least one positive entry.
For such multi-type GWI processes, under second-order moment conditions on the offspring and immigration distributions,
 they showed the weak convergence of the appropriately scaled one-dimensional distributions.
In particular, in the case that we consider in our present paper, namely, when a critical decomposable $p$-type
 GWI process $(\bX_k)_{k\geq 0} = ([X_{k,1},\dots,X_{k,p}]^\top)_{k\geq 0}$
has a lower triangular offspring mean matrix $(a_{i,j})_{i,j=1}^p\in[0,\infty)^{p \times p}$,
where $a_{i,i}=1$ for each $i\in\{1,\dots,p\}$ and $a_{i+1,i}>0$ for each $i\in\{1,\dots,p-1\}$,
Foster and Ney \cite[Theorems 4 and 5]{FosNey2} showed that $([n^{-1} X_{n,1}, \dots, n^{-p}X_{n,p}]^\top)_{n>0}$
converges in distribution as $n\to\infty$, and they characterized the limit distribution by its Laplace transform.
In our Theorem \ref{main_1}, we extend this result of Foster and Ney \cite{FosNey2} by proving weak convergence of
a sequence of appropriately scaled random step processes formed from $(\bX_k)_{k\geq 0}$,
and we describe the limit process using independent squared Bessel processes and 
 some linear combinations of their (at most $(p-1)$-fold) iterated integral processes.
We extend the results of Foster and Ney \cite[Section 9]{FosNey2} in another respect as well, namely,
we allow that the entries $a_{i+1,i}$, $i\in\{1,\dots,p-1\}$
of the offspring mean matrix $(a_{i,j})_{i,j=1}^p$ be $0$ as well, and in all the possible cases, we prove not only weak
convergence of one-dimensional distributions, but weak convergence of a sequence of appropriate random step processes
formed from the branching processes in question as well, see Theorem \ref{main_1}.
For a detailed comparison of the results of Foster and Ney \cite[Theorems 4 and 5]{FosNey2} and our Theorem \ref{main_1}, see Remark \ref{Rem_FosNey}.

In Barczy et al.\ \cite{BarBezPap2}, we described the asymptotic behavior of a critical decomposable 2-type GWI process.
Under second- or fourth-order moment assumptions on the offspring and immigration distributions, we proved that
 a sequence of appropriately scaled random step processes formed from a critical decomposable
 2-type GWI process converges weakly.
The limit process can be described using one or two independent squared Bessel processes and possibly the unique stationary
 distribution of an appropriate single-type subcritical GWI process.
Our results completed and extended the results of Foster and Ney \cite{FosNey2} for some strongly critical decomposable 2-type GWI process.

In Barczy and Bezd\'any \cite{BarBez}, we described the asymptotic behavior of a critical decomposable 3-type GWI process,
whose offspring mean matrix is triangular with all diagonal entries $1$.
Under second- or fourth-order moment assumptions on the offspring and immigration distributions, we proved that
a sequence of appropriately scaled random step processes formed from such a GWI process converges weakly.
The coordinate processes of the limit process are given as independent squared Bessel processes
 and some linear combinations of their 1-fold and 2-fold iterated integral processes.

For a detailed comparison of the result of this paper and those of Barczy et al. \cite{BarBezPap2} and Barczy and Bezd\'any \cite{BarBez}, see Remark \ref{Rem_spec}.

Concerning limit theorems for critical decomposable multi-type Galton-Watson processes, recently,
 Sagitov \cite[Theorem 4]{Sag} has proved some conditional limit theorems.
Namely, he has considered a sequence of critical decomposable $p$-type Galton--Watson processes (without immigration)
 with a common offspring mean matrix $(a_{i,j})_{i,j=1}^p$ that is upper triangular such that $a_{i,i}=1$, $i=1,\ldots,p$,
 and $a_{i-1,i}>0$, $i=2,\ldots,p$, and assumed that the $n^{\mathrm{th}}$ process starts with $n$ individuals of type 1.
He has shown weak convergence of the finite dimensional distributions of appropriately scaled random step processes formed
 from the given Galton-Watson processes as $n\to\infty$.
The limit distributions in Theorem 4 in Sagitov \cite{Sag} can be reformulated as the finite dimensional distributions
 of a $p$-dimensional stochastic process of which the $i^{\mathrm th}$ coordinate process
 is the $(i-1)$-th iterated integral process of its first coordinate process,
 multiplied by some appropriate constants, where $i\in\{2,\dots,p\}$.
In certain cases, the limit process in our Theorem \ref{main_1} has a very similar structure, see Remarks \ref{Rem_altforms} and \ref{Rem_FosNey} for more details.
We also mention that a similar type of limit distribution as in Theorem 4 in Sagitov \cite{Sag} appears 
 also in a limit theorem for critical Galton-Watson processes with
 overlapping generations (such branching processes can be regarded as discrete time versions of the Crump-Mode-Jagers process)
  with an individual life time distribution $L$ such that $\NN\ni n\mapsto\sum_{j=0}^n \PP(L>j)$ is regularly varying at infinity with exponent in $[0,1]$,
 see Sagitov \cite[Corollary 3 and remark (2) on page 89]{Sag}.

Decomposable multi-type Galton–Watson processes, 
 with or without immigration, have found some applications as well.
They arise naturally as stochastic models for structured populations,
 in particular, one can use them for modelling 
 the reproduction of biological populations consisting 
 of two types of individuals which display an
 irreversibility property in the sense that individuals 
 of one type might give birth to descendants of both kinds, 
 whereas those of the other type can have descendants only of their own kind
 (see Jagers \cite{Jag}).
One can use them to analyze the interactions between division,
 mutation and selection in evolutionary models 
 (see Balelli et al.\ \cite{BalMilWai})
 or to model the outcomes of polymerase chain reaction (PCR) bar encodings 
 (see Sagitov and St{\aa}hlberg \cite{SagStaahl}).
For a more detailed exposition of possible applications,
 see the Introduction of Barczy and Bezd\'any \cite{BarBez}.

The paper is organized as follows.
In Section \ref{Section_multi-type_GWI}, we introduce some notations, as well as we recall the definition
and some basic properties (e.g. mean function, classification and decomposability) of multi-type GWI processes.
In Section \ref{Section_conv_results}, we specify the class of GWI processes that this paper focuses on,
namely, we consider critical multi-type GWI processes whose offspring mean matrix is lower triangular with all diagonal entries $1$.
Section \ref{Section_conv_results} contains also our main result, Theorem \ref{main_1},
 a functional limit theorem for a critical decomposable $p$-type GWI process in question.
Under second-order moment conditions on the offspring and immigration distributions, we show that
a sequence of appropriately scaled random step processes formed from the GWI process in question converges weakly.
We describe the limit process using independent squared Bessel processes (given as pathwise unique strong solutions of stochastic differential equations (SDEs)) 
 and some linear combinations of their (at most $(p-1)$-fold) iterated integral processes.
We also give alternative representations of the limit process (see Remark \ref{Rem_altforms}),
 and compare our results with those of Wei and Winnicki \cite[Theorem 2.1]{WW}, Barczy et al. \cite{BarBezPap2}, Barczy and Bezd\'any \cite{BarBez},
 and Foster and Ney \cite{FosNey2} (see Remarks \ref{Rem_spec} and \ref{Rem_FosNey}).
Section \ref{Proof1} contains the proof of Theorem \ref{main_1}, which we split into three main steps.
The proof relies on functional limit theorems for sequences of random step processes formed from
 martingale differences, careful applications of corollaries
 of the continuous mapping theorem (see Lemmas \ref{Conv2Funct}, \ref{time_transform} and \ref{ConvIteratedIntegral}),
 and estimates on the (mixed) moments of the coordinates of multi-type GWI processes and some related martingale differences.
The paper ends with four appendices.
Appendix \ref{appInt} contains a result about sums of binomially weighted values of a function defined on the set of positive integers,
and some properties of certain triangular matrices (recalled from Barczy and Bezd\'any \cite{BarBez}).
In Appendix \ref{GWI_moments}, we recall some formulae and estimates for (mixed) moments of the coordinates of multi-type GWI processes
from Isp\'any and Pap \cite{IspPap2} and our recent paper Barczy and Bezd\'any \cite{BarBez},
and further, we also prove a new estimate (see Lemma \ref{M_sum_growth}).
In Appendix \ref{CMT}, we recall versions of the continuous mapping theorem (see Lemmas \ref{lemma:kallenberg} and \ref{Conv2Funct}) and present some
corollaries of it (see Lemmas \ref{time_transform} and \ref{ConvIteratedIntegral}) that will be useful in the proof of Theorem \ref{main_1}.
Appendix \ref{section_conv_step_processes} contains some results on the convergence of random step processes and critical single-type GWI processes,
due to Isp\'any and Pap \cite{IspPap} and Wei and Winnicki \cite[Theorem 2.1]{WW}, respectively.

The results of this paper generalize and unify those of Wei and Winnicki \cite{WW}, Barczy et al.\ \cite{BarBezPap2} and Barczy and Bezd\'any \cite{BarBez}
 concerning critical single-type, $2$-type and $3$-type decomposable GWI processes, respectively (for more details, see Remark \ref{Rem_spec}).
Our results also complete and extend some results of Foster and Ney \cite{FosNey2} (for more details, see Remark \ref{Rem_FosNey}).
Compared to our previous results in Barczy et al.\ \cite{BarBezPap2} and Barczy and Bezd\'any \cite{BarBez} for $2$-type and $3$-type GWI processes, respectively,
 in our present paper, we have found more compact formulae for the limit process in Theorem \ref{main_1},
 we have reduced the previous fourth-order moment assumptions on the offspring and immigration distributions to second-order ones,
 and the proof contained in Section \ref{Proof1} is more concise than the one that can be found in Barczy and Bezd\'any \cite{BarBez}.

\section{Multi-type GWI processes}\label{Section_multi-type_GWI}

In this section, first, we collect some notations that are used throughout the paper, and then we recall the definition, some moment formulae, and 
 the notions of decomposability and strong criticality for multi-type GWI processes.

Let $\ZZ_+$, $\NN$, $\RR$, $\RR_+$ and $\RR_{++}$ denote the set
 of non-negative integers, positive integers, real numbers, non-negative real
 numbers and positive real numbers, respectively.
For $x,y\in\RR$, the minimum and maximum of $x$ and $y$ are denoted by $x\wedge y$  and $x\vee y$, respectively.
For $k,m\in\ZZ_+$ with $k<m$, we define $\binom{k}{m}:=0$.
For $i,j\in\ZZ_+$, let $\delta_{i,j}:=0$ if $i\ne j$, and $\delta_{i,j}:=1$ if $i=j$ (known as the Kronecker delta).
 For functions $f:\NN\to\RR$ and $g:\NN\to\RR$, the notation $f(k) = \OO(g(k))$, $k\in\NN$, means that
 there exists a positive real number $M$ such that $\vert f(k)\vert\leq M g(k)$ for each $k\in\NN$.
More generally, for $d\in\ZZ_+$ and functions $G:\RR^{d+1}\to\RR$, $g_1,\dots,g_d:\NN\to\RR$, and $h:\NN\to\RR$, the notation
 \begin{equation}\label{help_bigO}
  G(k,\OO(g_1(k)),\dots,\OO(g_d(k)))=\OO(h(k)),\qquad k\in\NN,
 \end{equation}
means that for all $c_1,\dots,c_d\in\RR_+$,
 there exists a positive real number $M$ such that
\begin{equation*}
|G(k, c_1 g_1(k), \dots, c_d g_d(k))|\leq M h(k),\qquad k\in\NN.
\end{equation*}
For example, $\OO(1)=\OO(k)$, $k\in\NN$, because for all $c\in\RR_+$, $|c\cdot1|\leq c\cdot k$, $k\in\NN$,
but we point out that $\OO(k)=\OO(1)$, $k\in\NN$, does not hold.
The relation $\OO(1)=\OO(k)$, $k\in\NN$, fits into the framework \eqref{help_bigO} with the following choices:
 $d=1$, $G:\RR^2\to \RR$, $G(x,y)=y$, $(x,y)\in\RR^2$, $g_1:\NN\to\RR$, $g_1(k)=1$, $k\in\NN$, and
 $h:\NN\to\RR$, $h(k)=k$, $k\in\NN$.
The relation $f(k) = \OO(g(k))$, $k\in\NN$, also fits into the framework \eqref{help_bigO} with the following choices:
 $d=1$, $G:\RR^2\to \RR$ satisfying $G(k,y)=f(k)$, $(k,y)\in\NN\times\RR$, $g_1$ can be arbitrary,
 and $h:\NN\to\RR$, $h(k)=g(k)$, $k\in\NN$; or with the choices
 $d=1$, $G:\RR^2\to \RR$, $G(x,y)=y$, $(x,y)\in\RR^2$, $g_1:\NN\to\RR$, $g_1(k):=f(k)$, $k\in\NN$,
 and $h:\NN\to\RR$, $h(k)=g(k)$, $k\in\NN$.
For functions  $f:\NN\to\RR$ and $g:\NN\to\RR$, by the notation $f(k)=\OO(g(k))\to0$ as $k\to\infty$,
  we mean that $f(k)=\OO(g(k))$, $k\in\NN$, and $\lim_{k\to\infty}g(k)=0$, which imply that $\lim_{k\to\infty}f(k)=0$.

The Euclidean norm on $\RR^d$ is denoted by $\Vert\cdot\Vert$, where $d\in\NN$.
The $d\times d$ identity matrix is denoted by $\bI_d$.
For a matrix $\bB\in\RR^{d\times d}$ and $m\in\ZZ_+$, the $(i,j)$-th entry of $\bB^m$ is denoted by $b_{i,j}^{[m]}$,
 where $i,j\in\{1,\ldots,d\}$.
For a function $f:\RR\to\RR$, its positive part is denoted by $f^+$.
For $m\in\ZZ_+$ and a function $f: \cH\to\cH$, where $\cH$ is some set, we denote the
$m$-fold composition of $f$ with itself by $f^{\circ(m)}$, that is, $f^{\circ(0)}$ is the identity map on $\cH$,
 and $f^{\circ(m+1)}:=f\circ f^{\circ(m)}$, $m\in\ZZ_+$.
Every random variable will be defined on a fixed probability space $(\Omega, \cA, \PP)$.
Convergence in probability and equality in distribution
are denoted by $\stoch$ and $\distre$, respectively.
For $d\in\NN$, we will use $\DD(\RR_+,\RR^d)$ for the set of c\`adl\`ag functions from $\RR_+$ to $\RR^d$,
and $\CC(\RR_+,\RR^d)$ for the set of continuous functions from $\RR_+$ to $\RR^d$.
For $d\in\NN$, the convergence in the sense of the Skorokhod $J_1$
topology of functions in $\DD(\RR_+,\RR^d)$,
the weak convergence of the finite dimensional distributions of $\RR^d$-valued stochastic processes  with sample paths in $\DD(\RR_+, \RR^d)$,
and the weak convergence of $\RR^d$-valued stochastic processes with sample paths in $\DD(\RR_+, \RR^d)$
are denoted by $\skor{d}$, $\distrf$, and $\distr$, respectively
(for more details and notations, e.g., for the metric that generates the Skorokhod $J_1$ topology
and for the definition of locally uniform convergence $\lu$, see Appendix \ref{CMT}).

Next, we recall the definition and first-order moment formulae of $p$-type GWI processes, where $p\in\NN$.
For each $k \in \ZZ_+$ and $i \in \{ 1, \dots, p \}$, the number of
 individuals of type $i$ in the $k^\mathrm{th}$ generation is denoted by $X_{k,i}$.
For simplicity, we suppose that the initial values are $X_{0,i} = 0$, $i \in \{ 1, \dots, p \}$.
By $\xi_{k,j,i,\ell}$ we denote the number of type $\ell$ offspring
 produced by the $j^\mathrm{th}$ individual who is of type $i$
 belonging to the $(k-1)^\mathrm{th}$ generation.
The number of type $i$ immigrants in the $k^\mathrm{th}$ generation is denoted by $\vare_{k,i}$.
Consider the random vectors
 \[
   \bX_k := \begin{bmatrix}
             X_{k,1} \\
             \vdots \\
             X_{k,p}
            \end{bmatrix} , \qquad
   \bxi_{k,j,i} := \begin{bmatrix}
                    \xi_{k,j,i,1} \\
                    \vdots \\
                    \xi_{k,j,i,p}
                   \end{bmatrix} , \qquad
   \bvare_k := \begin{bmatrix}
                \vare_{k,1} \\
                \vdots \\
                \vare_{k,p}
               \end{bmatrix} .
 \]
Then we have
 \begin{equation}\label{MBPI(d)}
  \bX_k = \sum_{i=1}^p \sum_{j=1}^{X_{k-1,i}} \bxi_{k,j,i} + \bvare_k , \qquad
  k \in \NN,
 \end{equation}
 with $\bX_0=\bzero$ (and using the convention $\sum_{i=1}^0:=\bzero$).
Here $\big\{\bxi_{k,j,i}, \, \bvare_k : k, j \in \NN, \, i \in \{ 1, \dots, p \} \big\}$
  are supposed to be independent.
Moreover, $\big\{\bxi_{k,j,i} : k, j \in \NN\big\}$ for each
 $i \in \{1, \dots, p\}$, and $\{\bvare_k : k \in \NN\}$ are
 supposed to consist of identically distributed $\ZZ_+^p$-valued random vectors.
For notational convenience, let $\{\bxi_i : i \in \{1, \ldots, p\}\}$ and $\bvare$ be random vectors
 such that $\bxi_i \distre \bxi_{1,1,i}$ for each $i \in \{1, \ldots, p\}$ and $\bvare \distre \bvare_1$.

In all what follows, we will suppose
\begin{equation}\label{baseconds}
\EE(\|\bxi_i\|^2) < \infty, \quad i\in\{ 1,\dots, p\}, \quad\text{and}\quad \EE (\|\bvare\|^2) < \infty.
\end{equation}
Introduce the notations
 \begin{gather}\label{help_jelolesek}
  \begin{split}
 \bA := \begin{bmatrix}
          \EE(\bxi_1) & \cdots &           \EE(\bxi_p)
         \end{bmatrix} \in \RR^{p \times p}_+ , \qquad
  \bb := \EE(\bvare) \in \RR^p_+ , \\
  \bV^{(i)} := \var(\bxi_i) \in \RR^{p \times p} , \quad i\in\{1,\dots,p\},\qquad
  \bV^{(0)} := \var(\bvare) \in \RR^{p \times p}.
  \end{split}
 \end{gather}
The matrix $\bA$ and the vector $\bb$ are called the offspring mean matrix and the immigration mean vector, respectively.
Note that $\bA=(a_{i,j})_{i,j=1}^p$ with $a_{i,j} := \EE(\xi_{1,1,j,i})$, $i,j\in\{1,\ldots,p\}$,
 and we mention that some authors define the offspring mean matrix as \ $\bA^\top$.
Note also that $\bV^{(r)}=(v_{i,j}^{(r)})_{i,j=1}^p$ with $v_{i,j}^{(r)}:=\cov(\xi_{1,1,r,i},\xi_{1,1,r,j})$, $i,j,r\in\{1,\dots,p\}$,
 and $\bb=(b_i)_{i=1}^p$ with $b_i:=\EE(\vare_{1,i})$, $i\in\{1,\dots,p\}$.

For each $k \in \ZZ_+$, let $\cF_k^\bX := \sigma(\bX_0,\dots, \bX_k)$, where $\cF_0^\bX=\sigma(\bzero)=\{\emptyset,\Omega\}$.
By \eqref{MBPI(d)}, we get
 \begin{equation}\label{mart}
  \EE(\bX_k \mid \cF_{k-1}^\bX)
  = \sum_{i=1}^p X_{k-1,i} \EE(\bxi_i) + \bb
  = \bA \bX_{k-1} + \bb,\qquad k\in\NN .
 \end{equation}
Consequently,
 \begin{equation*}
  \EE(\bX_k) = \bA \EE(\bX_{k-1}) + \bb , \qquad k \in \NN ,
 \end{equation*}
and, since $\bX_0=\bzero$, we have
 \begin{equation*}
  \EE(\bX_k) = \sum_{j=0}^{k-1} \bA^j \bb , \qquad k \in \ZZ_+.
 \end{equation*}
The above equation shows that the offspring mean matrix $\bA$ plays a key role
 in the behavior of the expected values of the process $(\bX_k)_{k\in\ZZ_+}$.
A $p$-type GWI process $(\bX_k)_{k\in\ZZ_+}$ is referred to respectively as \emph{subcritical}, \emph{critical} or
 \emph{supercritical} if $\varrho(\bA) < 1$, $\varrho(\bA) = 1$
 or $\varrho(\bA) > 1$, where $\varrho(\bA)$ denotes the spectral radius of the matrix $\bA$
 (see, e.g., Athreya and Ney \cite[V.3]{AN} or Quine \cite{Q}), and it is called indecomposable or decomposable
 if $\bA$ is irreducible or reducible, respectively.
Recall that $\bA$ is called reducible if there exists a permutation matrix $\bP\in\RR^{p\times p}$
 and $q\in\{1,\ldots,p-1\}$ (and hence $p$ should be at least $2$) such that
\begin{equation}\label{reducMatrix}
\bA=\bP\begin{bmatrix}\bR&\bzero\\ \bS & \bT\end{bmatrix}\bP^\top,
\end{equation}
where $\bR\in\RR^{q\times q},\bS\in\RR^{(p-q)\times q},\bT\in\RR^{(p-q)\times(p-q)}$, and $\bzero\in\RR^{q\times(p-q)}$ is a null matrix.
The matrix $\bA$ is called irreducible if it is not reducible, see, e.g., Horn and Johnson \cite[Definitions 6.2.21 and 6.2.22]{HJ}.
We do emphasize that no 1-by-1 matrix is reducible.
If $p\in\NN$ and $\bA\in\RR_+^{p\times p}$, then there exist
a permutation matrix $\bP\in\RR_+^{p\times p}$ and $N\in\{1,\dots,p\}$ such that
\begin{equation}\label{reducible}
\bA=\bP\begin{bmatrix}
\bA_{1,1} & \bzero & \dots & \bzero\\
\bA_{2,1} & \bA_{2,2} & \dots & \bzero\\
\vdots & \vdots & \ddots & \vdots\\
\bA_{N,1} & \bA_{N,2} & \dots & \bA_{N,N}
\end{bmatrix}\bP^\top,
\end{equation}
where $\bA_{i,j}\in\RR^{p_i\times p_j}$ for some $p_i,p_j\in\{1,\dots,p\}$, $i,j\in\{1,\dots,N\}$, $p_1+\dots+p_N=p$,
and $\bA_{i,i}$, $i\in\{1,\dots,N\}$ are all irreducible, see Barczy and Bezd\'any \cite[Section 2]{BarBez}.
In particular, if $\bA$ is irreducible, then $\bP=\bI_p$ and $N=1$.
We say that $(\bX_k)_{k\in\ZZ_+}$ is strongly critical if $\varrho(\bA_{i,i})=1$ for each $i\in\{1,\ldots,N\}$,
 where $\bA_{i,i}$, $i\in\{1,\ldots,N\}$, are the irreducible matrices in the decomposition \eqref{reducible} of $\bA$.

\section{Asymptotic behavior of random step processes formed from
  some strongly critical decomposable $p$-type GWI processes}
\label{Section_conv_results}

In what follows, let $(\bX_k)_{k\in\ZZ_+}$ be a $p$-type Galton-Watson process with
 immigration starting from $\bX_0=\bzero$ such that the moment condition \eqref{baseconds} holds
 and that the offspring mean matrix $\bA$ takes the form
\begin{equation}\label{strongcrit_triang}
\bA=\begin{bmatrix}
1 & 0 & 0 & \dots & 0\\
a_{2,1} & 1 & 0 & \dots & 0\\
a_{3,1} & a_{3,2} & 1 & \dots & 0\\
\vdots & \vdots & \vdots & \ddots & \vdots\\
a_{p,1} & a_{p,2} & a_{p,3} & \dots & 1
\end{bmatrix}\in\RR_+^{p\times p}.
\end{equation}
Clearly, $\varrho(\bA)=1$, thus $(\bX_k)_{k\in\ZZ_+}$ is critical, and if $p\geq2$, then it is decomposable.
Furthermore, if $p\geq2$, then $(\bX_k)_{k\in\ZZ_+}$ is a special strongly critical decomposable $p$-type GWI process.
Indeed, $\bA$ is of the form \eqref{reducible} with $\bP=\bI_p$, $N=p$, $p_1=\dots=p_p=1$, and $\bA_{i,i}:=[1]$.
Note also that for each $i\in\{1,\dots,p\}$, if $a_{i,j}=\delta_{i,j}$, $j\in\{1,\dots,p\}$,
 then $(X_{k,i})_{k\in\ZZ_+}$ is a critical single-type GWI process, since $a_{i,j}=0$, $j\in\{1,\ldots,p\}\setminus\{i\}$
 implies that $\PP(\xi_{1,1,j,i}=0)=1$, $j\in\{1,\ldots,p\}\setminus\{i\}$, and hence
 $X_{k,i} = \sum_{j=1}^{X_{k-1,i}} \xi_{k,j,i,i} + \vare_{k,i}$, $k\in\NN$, with $X_{0,i}=0$ and $\EE(\xi_{1,1,i,i})=1$.
In particular, the first coordinate process $(X_{k,1})_{k\in\ZZ_+}$ of $(\bX_k)_{k\in\ZZ_+}$ is a critical (single-type) GWI process.

For a matrix $\bA\in\RR_+^{p\times p}$ which satisfies \eqref{strongcrit_triang},
 let us introduce the notation
 \begin{equation}\label{etadef}
 \eta_i(\bA) :=\max_{m\in\{1,\dots,p\}}\{m \mid \exists\, j\in\{1,\dots,p\}: c_{i,j}^{[m-1]}>0\},
   \qquad i\in\{1,\ldots,p\},
 \end{equation}
 where $c_{i,j}^{[m]}$, $i,j\in\{1,\dots,p\}$, $m\in\ZZ_+$, denotes the $(i,j)$-th entry of $\bC^{m}$ with $\bC:=\bA-\bI_p$.
When it is clear from the context which matrix $\bA$ is involved, we will write $\eta_i:=\eta_i(\bA)$, $i\in\{1,\dots,p\}$.
Due to $c_{i,i}^{[0]}=(\bI_p)_{i,i}=1$, we have $\eta_i(\bA)\geq1$ for all $\bA$ satisfying \eqref{strongcrit_triang} and each $i\in\{1,\dots,p\}$.

We check that, without loss of generality, one can suppose the existence of a $q\in\{1,\dots,p\}$ such that for each $i\in\{1,\dots,p\}$,
 we have $\eta_i=1$ if and only if $i\leq q$, where $\eta_i$ is defined in \eqref{etadef}.
Note that this means that the upper-left $q\times q$ submatrix of $\bA$ is $\bI_q$.
More precisely, we check that if $\bA$ takes the form \eqref{strongcrit_triang}, then we can find an
 appropriate $p\times p$ permutation matrix $\bP$ such that the mean matrix $\bP\bA\bP^\top$ of
 the GWI process $(\bP\bX_k)_{k\in\ZZ_+}$ has the form \eqref{strongcrit_triang} and
 there exists a $q\in\{1,\dots,p\}$ such that for each $i\in\{1,\dots,p\}$,
 we have $\eta_i(\bP\bA\bP^\top)=1$ if and only if $i\leq q$.
First,
 note that $\eta_1(\bA)=1$, since $\bC=\bA-\bI_p$ is a lower triangular matrix such that all the entries
 in its first row are $0$, yielding that $c_{1,j}^{[\ell]}=0$, $j\in\{1,\ldots,p\}$, $\ell\in\NN$.
Note also that, for some $i\in\{2,\ldots,p\}$, $\eta_i=1$ if and only if $a_{i,j}=0$, $j\in\{1,\ldots,i-1\}$.
Let
 \[
    q_0:=\max\Big\{i\in\{1,\dots,p\} : \eta_r=1,\; \forall \,r\in\{1,\dots,i\}\Big\},
 \]
  which always exists due to $\eta_1=1$.
If $q_0=p$, then we have $q:=q_0=p$.
If $q_0<p$ and $\eta_r>1$ for each $r\in\{q_0+1,\dots,p\}$, then $q:=q_0$.
Otherwise, we have $q_0<p$ with $\eta_{q_0+1}>1$ ($\eta_{q_0+1}=1$ cannot hold,
 since it would lead to a contradiction with the definition of $q_0$),
 and there exists a smallest $r\in\{q_0+2,\dots,p\}$ such that $\eta_r=1$.
Now consider the $p\times p$ permutation matrix $\bP_1$ such that
 \begin{equation*}
 \bP_1\bX_k
   =\begin{bmatrix}X_{k,1} & \dots & X_{k,q_0} & X_{k,r} & X_{k,q_0+1} & \dots & X_{k,r-1} & X_{k,r+1} & \dots & X_{k,p}\end{bmatrix}^\top
 \end{equation*}
 for each $k\in\ZZ_+$.
 It is easy to see that $(\bP_1\bX_k)_{k\in\ZZ_+}$ is a $p$-type GWI process with an offspring mean matrix $\bA_1:=\bP_1\bA\bP_1^\top$,
 and we have that $\eta_i(\bA_1)=1$ for each $i\in\{1,\dots,q_0+1\}$.
We can now define $q_1$ for $\bA_1$ in a similar manner as $q_0$ was defined for $\bA$, and get that $q_1\geq q_0+1>q_0$.
If this value $q_1$ is not a suitable choice for $q$, then we can repeat this procedure at most $p-1$ times
 until all of the initial types $i\in\{1,\dots,p\}$ with $\eta_i(\bA)=1$ are permuted to the first $q$ positions,
and we have $q:=q_r\leq p$ for some $r\in\{0,\dots,p-1\}$.

 Recall that the entries of the matrices $\bV^{(1)},\ldots,\bV^{(p)}$ (introduced in \eqref{help_jelolesek})
 are denoted by $v_{i,j}^{(1)},\ldots, v_{i,j}^{(p)}$, $i,j\in\{1,\ldots,p\}$.
 Further, the coordinates of $\bb = \EE(\bvare)$ are denoted by $b_1,\ldots, b_p$.

\begin{Thm}\label{main_1}
Let $(\bX_k)_{k\in\ZZ_+}$ be a $p$-type GWI process such that
its offspring mean matrix $\bA$ satisfies \eqref{strongcrit_triang},
 $\EE(\|\bvare\|^2)<\infty$,  $\EE(\|\bxi_i\|^2)<\infty$ for each $i\in\{1,\dots,p\}$,
 and let  $q\in\{1,\ldots,p\}$ be such that, for each $j\in\{1,\ldots,p\}$, we have that $\eta_j=1$ if and only if $j\leq q$,
with $\eta_i:=\eta_i(\bA)$, $i\in\{1,\dots,p\}$, defined in \eqref{etadef}.
Then
 \begin{equation}\label{Conv_X_1}
  \left(\begin{bmatrix}
          n^{-\eta_1} X_{\nt,1} \\
          \vdots\\
          n^{-\eta_p} X_{\nt,p}
         \end{bmatrix}\right)_{t\in\RR_+}
  \distr
  \left(\begin{bmatrix}
          \cX_{t,1} \\
          \vdots\\
          \cX_{t,p}
         \end{bmatrix}\right)_{t\in\RR_+}
  \qquad \text{as \ $n \to \infty$,}
 \end{equation}
 where the limit process  $([\cX_{t,1}, \dots, \cX_{t,p}]^\top)_{t\in\RR_+}$ is given as follows:

\begin{itemize}
\item[(a)] The process $([\cX_{t,1}, \dots, \cX_{t,q}]^\top)_{t\in\RR_+}$ is the pathwise unique strong solution of the SDE
 \begin{equation}\label{SDE_X_1}
  \begin{cases}
   &\dd\cX_{t,1}
   = b_1 \, \dd t
     + \sqrt{v^{(1)}_{1,1} \, \cX_{t,1}^+} \,
       \dd \cW_{t,1} , \\[1mm]
   &\phantom{\dd\cX_{t,1}\;\;}\vdots\\[1mm]
   &\dd\cX_{t,q}
   = b_{q} \, \dd t
     + \sqrt{v^{(q)}_{q,q} \, \cX_{t,q}^+} \,
       \dd \cW_{t,q} ,
  \end{cases}
  \qquad t \in\RR_+ ,
 \end{equation}
 with initial value $[\cX_{0,1} , \dots, \cX_{0,q}]^\top = \bzero\in\RR^q$, where $(\cW_{t,1})_{t\in\RR_+},\dots,(\cW_{t,q})_{t\in\RR_+}$
 are independent standard Wiener processes yielding the independence of $(\cX_{t,1})_{t\in\RR_+},\dots, (\cX_{t,q})_{t\in\RR_+}$ as well.

\item[(b)] For each $i\in\{q+1,\dots,p\}$ and $m\in\{1,\dots,\eta_i-1\}$, we have that
 \begin{align}\label{SDE_X_2}
  \begin{split}
 \cX_{t,i} &=\sum_{j=1}^qc_{i,j}^{[\eta_i-1]}\int_0^t\left(\int_0^{s_{\eta_i-1}}\dots
  \left( \int_0^{s_3}\left(\int_0^{s_2}\cX_{s_1,j}\,\dd s_1\right)\dd s_2\right) \dots \, \dd s_{\eta_i-2}\right)\,\dd s_{\eta_i-1}\\
	& =\sum_{ \substack{j\in\{1,\dots,p\},\\\eta_j=\eta_i-m}}c_{i,j}^{[m]}\int_0^t\left(\int_0^{s_m}\dots
  \left( \int_0^{s_3}\left(\int_0^{s_2}\cX_{s_1,j}\,\dd s_1\right)\dd s_2\right) \dots\,\dd s_{m-1}\right)\,\dd s_m
  \end{split}
 \end{align}
 \end{itemize}
for all $t\in\RR_+$, where we recall that $c_{i,j}^{[\eta_i-1]}$ denotes the $(i,j)$-th entry of $\bC^{\eta_i-1} = (\bA-\bI_p)^{\eta_i-1}$.
\end{Thm}

In the next remark, we rewrite the stochastic processes $(\cX_{t,i})_{t\in\RR_+}$, $i\in\{q+1,\dots,p\}$, given in \eqref{SDE_X_2}
 into several other forms.

\begin{Rem}\label{Rem_altforms}
 Using the notation $\varphi$ (see Definition \ref{maps}), in Theorem \ref{main_1}, for each $i\in\{q+1,\dots,p\}$ and $m\in\{1,\dots,\eta_i-1\}$,
 we have that
\begin{align*}
  \cX_{t,i} & = \sum_{ \substack{j\in\{1,\dots,p\},\\\eta_j=\eta_i-m}}c_{i,j}^{[m]}\left(\varphi^{\circ(m)}\left(\left(\cX_{s,j}\right)_{s\in\RR_+}\right)\right)(t)
            & = \sum_{ \substack{j\in\{1,\dots,i-1\},\\\eta_j=\eta_i-m}}c_{i,j}^{[m]}\left(\varphi^{\circ(m)}\left(\left(\cX_{s,j}\right)_{s\in\RR_+}\right)\right)(t)
\end{align*}
for all $t\in\RR_+$, where the second equality follows from the fact that $c_{i,j}^{[m]} = 0$ for $j\geq i$ (see \eqref{C_0}).
Furthermore, by choosing $m=1$, for each $i\in\{q+1,\dots,p\}$, we also get
\begin{equation*}
 \cX_{t,i} =\sum_{\substack{j\in\{1,\dots,p\},\\\eta_j=\eta_i-1 }}a_{i,j}\int_0^t\cX_{s,j}\,\dd s
          = \sum_{j=1}^{i-1}\delta_{\eta_j,\eta_i-1}a_{i,j}\int_0^t\cX_{s,j}\,\dd s,\qquad t\in\RR_+,
\end{equation*}
 since for $i,j\in\{1,\dots,p\}$, $\eta_j=\eta_i-1$ implies $i\neq j$, and consequently $c_{i,j}^{[1]}=a_{i,j}$.
Moreover, the iterated integral processes appearing in \eqref{SDE_X_2} may be written as integrals of certain stochastic processes.
Using again the notation $\varphi$, for all $j\in\{1,\dots,p\}$, $k\in\NN$, and $t\in\RR_+$, we check that
\begin{align}\label{help_08}
 \begin{split}
 \left(\varphi^{\circ(k)}\left((\cX_{s,j})_{s\in\RR_+}\right)\right)(t)
		=\frac{1}{(k-1)!}\int_0^t(t-s)^{k-1}\cX_{s,j}\,\dd s
		=\frac{1}{k!}\int_0^t\cX_{t-s,j}\,\dd s^{k}.
 \end{split}
 \end{align}
The first equality in \eqref{help_08} follows (pathwise) from the Cauchy formula for repeated integration.
The second equality in \eqref{help_08} can be seen by substitution (set $u=t-s$) and
rewriting the integral as a Lebesgue-Stieltjes integral (note that $\dd u^k=k u^{k-1}\,\dd u$):
\begin{align*}
    \frac{1}{(k-1)!}\int_0^t(t-s)^{k-1}\cX_{s,j}\,\dd s
     = \frac{1}{(k-1)!} \int_0^t u^{k-1} \cX_{t-u,j}\,\dd u
     = \frac{1}{k!} \int_0^t \cX_{t-u,j}\,\dd u^k.
   \end{align*}
Finally, for each $i\in\{q+1,\ldots,p\}$, taking into account that $\eta_i\geq 2$, the formulae in \eqref{help_08}
can be used to rewrite the iterated integrals in the description of $(\cX_{t,i})_{t\in\RR_+}$ in \eqref{SDE_X_2} by choosing $k=m$,
where $m\in\{1,\dots,\eta_i-1\}$.
\proofend
\end{Rem}

In the next remark, we specialize Theorem \ref{main_1} to the cases $p=1,2$, and $3$,
 which were already considered in Wei and Winnicki \cite{WW}, Barczy et al.\ \cite{BarBezPap2},
 and Barczy and Bezd\'any \cite{BarBez}, respectively.

\begin{Rem}\label{Rem_spec}
In the case $p=1$, it is clear that we have $q=1$, and the statement of Theorem \ref{main_1}
 coincides with the statement of Theorem 2.1 in Wei and Winnicki \cite{WW} in the special case $X_0=0$.
Theorems 2.1 and 2.2 of Barczy et al.\ \cite{BarBezPap2} and the results of Barczy and Bezd\'any \cite{BarBez}
 can be retrieved as special cases of Theorem \ref{main_1}
 by choosing $p$ and $q$ according to the following table:
\begin{align*}
 \begin{tabular}{|c|c|c|c|}
  \hline
   & $p$ & $q$ \\
  \hline
  \cite[Theorem 2.1]{BarBezPap2} & $2$ & $2$ \\
  \hline
  \cite[Theorem 2.2]{BarBezPap2} & $2$ & $1$\\
  \hline
  \cite[Theorem 3.1]{BarBez} & $3$ & $3$\\
  \hline
  \cite[Theorem 3.2]{BarBez} & $3$ & $2$\\
  \hline
  \cite[Theorems 3.3 and 3.4]{BarBez} & $3$ & $1$\\
  \hline
 \end{tabular}\;.
\end{align*}
This is clear in the cases of Theorem 2.1 in \cite{BarBezPap2} and Theorem 3.1 in \cite{BarBez}, where
 the offspring mean matrix is $\bI_2$ and $\bI_3$, respectively.
In the case of Theorem 2.2 in \cite{BarBezPap2},
it is easy to see that we have
\begin{equation*}
\bC= \bA-\bI_2 = \begin{bmatrix}
0 & 0\\
a_{2,1} & 0
\end{bmatrix},
\end{equation*}
$\bC^0=\bI_2$, and $\bC^m=\bzero\in\RR^{2\times2}$ for each $m\in\NN\setminus\{1\}$.
This implies that $\eta_1=1$ and $\eta_2=2$,
and $c_{2,1}^{[\eta_2-1]}=c_{2,1}^{[1]}=a_{2,1}$.
Thus, in this case, Theorem \ref{main_1} yields that
\begin{equation*}
\left(\begin{bmatrix}
n^{-1}X_{\nt,1}\\
n^{-2}X_{\nt,2}\\
\end{bmatrix}\right)_{t\in\RR_+}
\distr
\left(
\begin{bmatrix}
\cX_{t,1}\\
a_{2,1}\int_0^t\cX_{s,1}\,\dd s\\
\end{bmatrix}
\right)_{t\in\RR_+} \qquad \text{as $n\to\infty$,}
\end{equation*}
where $(\cX_{t,1})_{t\in\RR_+}$ satisfies the SDE \eqref{SDE_X_1} (with $q=1$),
which is the same limit process as given by Theorem 2.2 in \cite{BarBezPap2}.
One may proceed similarly for the remaining 3-type cases, where the offspring mean matrix is not $\bI_3$.
In these cases, by (4.5) in the proof of Corollary 4.1 in Barczy and Bezd\'any \cite{BarBez}, we have
\begin{equation*}
\bC^0=\bI_3,
\qquad
\bC^1=\begin{bmatrix}
0 & 0 & 0\\
a_{2,1} & 0 & 0\\
a_{3,1} & a_{3,2} & 0
\end{bmatrix},
\qquad
\bC^2=\begin{bmatrix}
0 & 0 & 0\\
0 & 0 & 0\\
a_{3,2}a_{2,1} & 0 & 0
\end{bmatrix},
\end{equation*}
and $\bC^m=\bzero\in\RR^{3\times3}$ for each $m\in\NN\setminus\{1,2\}$.
Consequently, by Lemma 4.3 in Barczy and Bezd\'any \cite{BarBez}, in these cases, $\eta_i$, $i\in\{1,2,3\}$,
 take the following values:
\begin{align}\label{Table2}
\centering
 \begin{tabular}{|c|c|c|c|}
  \hline
   & Theorem 3.2 of \cite{BarBez} & Theorem 3.3 of \cite{BarBez} & Theorem 3.4 of \cite{BarBez} \\
\hline
   $\eta_1$ & $1$ & $1$ & $1$ \\
  \hline
   $\eta_2$ & $1$ & $2$ & $2$ \\
  \hline
   $\eta_3$ & $2$ & $2$ & $3$\\
  \hline
 \end{tabular}\;.
\end{align}
It is now easy to check that the statement of Theorem \ref{main_1} agrees with the results of Barczy and Bezd\'any \cite{BarBez}.
For example, by Theorem \ref{main_1} and \eqref{Table2}, in the case of Theorem 3.4 of \cite{BarBez}, we have
\begin{align*}
\left(\begin{bmatrix}
n^{-1}X_{\nt,1}\\
n^{-2}X_{\nt,2}\\
n^{-3}X_{\nt,3}
\end{bmatrix}\right)_{t\in\RR_+}
\distr
&\left(\begin{bmatrix}
\cX_{t,1}\\
c_{2,1}^{[1]}\int_0^t\cX_{s,1}\,\dd s\\
c_{3,1}^{[2]}\int_0^t\int_0^r\cX_{s,1}\,\dd s\,\dd r
\end{bmatrix}\right)_{t\in\RR_+}
=
\left(\begin{bmatrix}
\cX_{t,1}\\
a_{2,1}\int_0^t\cX_{s,1}\,\dd s\\
a_{3,2}a_{2,1}\int_0^t\int_0^r\cX_{s,1}\,\dd s\,\dd r
\end{bmatrix}\right)_{t\in\RR_+} \\
& =
\left(\begin{bmatrix}
\cX_{t,1}\\
a_{2,1}\int_0^t\cX_{s,1}\,\dd s\\
a_{3,2}\int_0^t\cX_{s,2}\,\dd s
\end{bmatrix}\right)_{t\in\RR_+}
\qquad \text{as $n\to\infty$,}
\end{align*}
 which indeed coincides with the statement of Theorem 3.4 of \cite{BarBez}.
We call attention to the fact that the second-order moment conditions on the offspring and immigration distributions 
 in Theorem \ref{main_1} are weaker than the corresponding fourth-order moment conditions
 in Theorem 2.1 of Barczy et al.\ \cite{BarBezPap2} and Theorems 3.1 and 3.2 of Barczy and Bezd\'any \cite{BarBez}.
\proofend
\end{Rem}

In the following remark we compare Theorem \ref{main_1} with some results of Foster and Ney \cite{FosNey2} for critical decomposable multi-type GWI processes.

\begin{Rem}\label{Rem_FosNey}
Theorem \ref{main_1} extends and completes some results of Foster and Ney \cite{FosNey2} for critical decomposable multi-type GWI processes.
In the special case when the offspring mean matrix $\bA$ satisfies \eqref{strongcrit_triang}
 with $a_{i+1,i}>0$, $i\in\{1,\dots,p-1\}$, Theorems 4 and 5 of Foster and Ney \cite{FosNey2} state that, under second-order moment conditions
 on the offspring and immigration distributions,  we have
\begin{equation}\label{FosNey_conv}
\begin{bmatrix}
          n^{-1} X_{n,1} \\
          \vdots\\
          n^{-p} X_{n,p}
         \end{bmatrix}
  \distr
\begin{bmatrix}
          \cX_{1} \\
          \vdots\\
          \cX_{p}
         \end{bmatrix}
  \qquad \text{as \ $n \to \infty$,}
 \end{equation}
where the distribution of $[\cX_1,\ldots,\cX_p]^\top$ is characterised by its Laplace transform.
We can derive this convergence as a consequence of Theorem \ref{main_1}.
Since $a_{i+1,i}>0$, $i\in\{1,\dots,p-1\}$, by Remark \ref{matrix_power_remark},
 we have that $c_{i,1}^{[i-1]}=\prod_{r=1}^{i-1} a_{r+1,r}>0$, $i\in\{2,\ldots,p\}$,
 and $c_{i,r}^{[i]}=0$ for each $i,r\in\{1,\dots,p\}$,
and consequently, using also that $c_{1,1}^{[0]}=1$, we get that $\eta_i=i$ for each $i\in\{1,\dots,p\}$, and $q=1$.
The convergence \eqref{FosNey_conv} is an immediate consequence of \eqref{Conv_X_1} (take $t=1$).
Moreover, we have that $\cX_i\distre\cX_{1,i}$ for $i\in\{1,\dots,p\}$, where $(\cX_{t,1})_{t\in\RR_+}$ is
a squared Bessel process (see \eqref{SDE_X_1}), and with the notation $\varphi$ (see Definition \ref{maps}), we have
\begin{equation*}
\cX_{t,i}= c_{i,1}^{[i-1]} \left(\varphi^{\circ(i-1)}\left(\left(\cX_{s,1}\right)_{s\in\RR_+}\right)\right)(t)
         = \left(\prod_{r=1}^{i-1} a_{r+1,r}\right) \left(\varphi^{\circ(i-1)}\left(\left(\cX_{s,1}\right)_{s\in\RR_+}\right)\right)(t)
\end{equation*}
 for $t\in\RR_+$ and $i\in\{2,\dots,p\}$.
It is worth noting that Theorems 4 and 5 of Foster and Ney \cite{FosNey2} are stated for a broader class of multi-type GWI processes
 than what we considered in this remark.
Our Theorem \ref{main_1} is also more general in the sense that we allow that $a_{i+1,i}=0$ for any $i\in\{1,\dots,p-1\}$,
 but is less general in the sense that we consider only multi-type GWI processes whose offspring mean matrix is triangular,
 instead of multi-type GWI processes whose offspring mean matrix is a particular kind of block triangular matrix
 (for more details, see Foster and Ney \cite[Introduction]{FosNey2} or the second paragraph of the present paper).
\proofend
\end{Rem}

\section{Proof of Theorem \ref{main_1}}\label{Proof1}%

We divide the proof into three main steps.

\smallskip

\textsl{Step 1}. A useful decomposition.

{\sl Step 1/(a).}
Let us introduce the sequence
 \begin{equation}\label{Mk1}
 \bM_k := \bX_k - \EE(\bX_k \mid \cF_{k-1}^{\bX})
   = \bX_k - \bA \bX_{k-1} - \bb, \qquad   k \in \NN ,
 \end{equation}
 of martingale differences with respect to the filtration
 \ $(\cF_k^\bX)_{k\in\ZZ_+}$, \ where we used \eqref{mart}.
From \eqref{Mk1}, we obtain the recursion
 $ \bX_k = \bA \bX_{k-1} + \bM_k + \bb$, $k \in \NN$,
 which together with \ $\bX_0=\bzero$ \ implies
 \begin{equation}\label{X}
   \bX_k = \sum_{\ell=1}^k \bA^{k-\ell} (\bM_\ell + \bb) , \qquad
   k \in \NN.
 \end{equation}

Consequently, by Lemma \ref{matrixlemma},
we have
 \begin{equation}\label{X2}
   \bX_k = \sum_{\ell=1}^k \sum_{m=0}^{p-1}\binom{k-\ell}{m}\bC^m (\bM_\ell + \bb)= \sum_{m=0}^{p-1}\bC^m\sum_{\ell=1}^k\binom{k-\ell}{m} (\bM_\ell + \bb), \qquad
   k \in \NN,
 \end{equation}
where recall that $\bC=\bA-\bI_p$.
Let us introduce the notation
\begin{equation}\label{M_sums}
X_{k,j}^{(m)}:=\sum_{\ell=1}^k\binom{k-\ell}{m}(M_{\ell,j}+b_j), \qquad m,k\in\ZZ_+, \quad j\in\{1,\dots,p\}.
\end{equation}
Note that $X_{0,j}^{(m)}=0$, $m\in\ZZ_+$, $j\in\{1,\ldots,p\}$.
Then \eqref{X2} can be written as
\begin{equation*}
   \bX_k=
\sum_{m=0}^{p-1}\bC^m\begin{bmatrix}
X_{k,1}^{(m)}\\
\vdots\\
X_{k,p}^{(m)}
\end{bmatrix}, \qquad k\in\NN.
\end{equation*}
Consequently, by the rules of matrix multiplication, for each $i\in\{1,\dots,p\}$, we have
\begin{equation}\label{decomp}
X_{k, i}=\sum_{m=0}^{p-1}\sum_{j=1}^pc_{i,j}^{[m]}X_{k,j}^{(m)}, \qquad k\in\NN,
\end{equation}
where recall that $c_{i,j}^{[m]}$ denotes the $(i,j)$-th entry of $\bC^m$, $i,j\in\{1,\dots,p\}$, $m\in\ZZ_+$.

{\sl Step 1/(b).}
We check that, for each $i\in\{1,\dots,p\}$, it holds that
 \begin{equation}\label{decomp2}
   X_{k, i} = \sum_{m=0}^{\eta_i-2}\sum_{j=1}^{p}c_{i,j}^{[m]}X_{k,j}^{(m)} + \sum_{j=1}^{q}c_{i,j}^{[\eta_i-1]}X_{k,j}^{(\eta_i-1)}, \qquad k\in\NN,
 \end{equation}
 where recall the convention $\sum_{\ell=0}^{-1}:=0$.
To show this, let $i\in\{1,\dots,p\}$ be fixed. By the definition of $\eta_i$ in \eqref{etadef},
 we have that $c_{i,j}^{[m]}=0$ for each $j\in\{1,\dots,p\}$ if $m > \eta_i-1$,
 thus all the terms in \eqref{decomp} corresponding to $m = \eta_i,\ldots,p-1$ are $0$, and hence
 we get that
 \begin{equation*}
  X_{k, i}=\sum_{m=0}^{\eta_i-1}\sum_{j=1}^pc_{i,j}^{[m]}X_{k,j}^{(m)}, \qquad k\in\NN.
 \end{equation*}
On the other hand, if $j\in\{q+1,\dots,p\}$ is such that $c_{i,j}^{[\eta_i-1]}>0$,
 then, by Lemma \ref{eta_ineq}, we have $\eta_i\geq\eta_j+\eta_i-1\geq2+\eta_i-1=\eta_i+1$,
 since due to $j>q$, we have $\eta_j>1$.
This leads us to a contradiction, so $c_{i,j}^{[\eta_i-1]}=0$ for $j\in\{q+1,\ldots,p\}$,
 and thus we obtain \eqref{decomp2} by
 \begin{align*}
 X_{k, i}=&\sum_{m=0}^{\eta_i-2}\sum_{j=1}^pc_{i,j}^{[m]}X_{k,j}^{(m)}+\sum_{j=1}^pc_{i,j}^{[\eta_i-1]}X_{k,j}^{(\eta_i-1)}\\
	    =&\sum_{m=0}^{\eta_i-2}\sum_{j=1}^pc_{i,j}^{[m]}X_{k,j}^{(m)}+\sum_{j=1}^qc_{i,j}^{[\eta_i-1]}X_{k,j}^{(\eta_i-1)}, \qquad k\in\NN.
 \end{align*}
Consequently, we may write
\begin{equation}\label{decomp_vector}
\bX_k=\begin{bmatrix}
X_{k,1}\\
\vdots\\
X_{k,p}
\end{bmatrix}=
\begin{bmatrix}
 \sum_{m=0}^{\eta_1-2}\sum_{j=1}^{p}c_{1,j}^{[m]}X_{k,j}^{(m)}\\
\vdots\\
 \sum_{m=0}^{\eta_p-2}\sum_{j=1}^{p}c_{p,j}^{[m]}X_{k,j}^{(m)}
\end{bmatrix}
+
\begin{bmatrix}
\sum_{j=1}^{q}c_{1,j}^{[\eta_1-1]}X_{k,j}^{(\eta_1-1)}\\
\vdots\\
\sum_{j=1}^{q}c_{p,j}^{[\eta_p-1]}X_{k,j}^{(\eta_p-1)}
\end{bmatrix}, \qquad k\in\NN.
\end{equation}
In Step 3, we will see that the second term in \eqref{decomp_vector}
 will serve as the “leading part” of $(\bX_k)_{k\in\ZZ_+}$,
 while the first term becomes negligible in the limit.

{\sl Step 1/(c).}
Note that, by \eqref{M_sums} and Lemma \ref{specSums}, for $j\in\{1,\dots,p\}$, $n\in\NN$ and $t\in\RR_+$, we have
 \begin{align}\label{integral_dec}
 \begin{split}
  X_{\nt,j}^{(m)}
             &= \sum_{\ell=1}^\nt\binom{\nt-\ell}{m}(M_{\ell,j}+b_j)
			 = n^m\left(\varphi_n^{\circ(m)}\left(\left(\sum_{\ell=1}^\ns(M_{\ell,j}+b_j)\right)_{s\in\RR_+}\right)\right)(t)\\
			 & = n^m\left(\varphi_n^{\circ(m)}\left(\left(X_{\ns,j}^{(0)}\right)_{s\in\RR_+}\right)\right)(t),
 \end{split}
 \end{align}
 where we used the notation $\varphi_n$ introduced in part (vi) of Definition \ref{maps}.
If $i\in\{1,\dots,q\}$, then, using that $c_{i,j}^{[m]}>0$ only if $j=i$ and $m=0$
 (following from the definitions of $q$ and $\eta_i$) and that $c_{i,i}^{[0]}=1$, formula \eqref{decomp} yields that
 $X_{k,i}=X_{k,i}^{(0)}$, $k\in\ZZ_+$, $i\in\{1,\ldots,q\}$.
Thus, by \eqref{integral_dec}, we get that
\begin{equation}\label{d_int}
X_{\nt,i}^{(m)}=n^m\left(\varphi_n^{\circ(m)}\left(\left(X_{\ns,i}\right)_{s\in\RR_+}\right)\right)(t), \qquad n\in\NN,\,m\in\ZZ_+,\, t\in\RR_+,\, i\in\{1,\ldots,q\}.
\end{equation}

\smallskip

\textsl{Step 2.} We describe the asymptotic behavior of the first $q$ coordinate processes of $(\bX_k)_{k\in\ZZ_+}$, namely,
we prove weak convergence of the processes $([n^{-1}X_{\nt,1}, \dots, n^{-1}X_{\nt,q}]^\top)_{t\in\RR_+}$, $ n\in\NN$.

{\sl Step 2/(a).}
We introduce the notations
\begin{gather}\label{q_notation}
\begin{split}
\bb_q:=\begin{bmatrix}
b_1\\
\vdots\\
b_q
\end{bmatrix},
\qquad
\bV_q^{(r)}:=\begin{bmatrix}
v_{1,1}^{(r)} & \dots & v_{1,q}^{(r)}\\
\vdots & \ddots & \vdots\\
v_{q,1}^{(r)} & \dots & v_{q,q}^{(r)}
\end{bmatrix}, \qquad r\in\{0,1,\dots,p\},\\
\bM_{k,q}:=\begin{bmatrix}
M_{k,1}\\
\vdots\\
M_{k,q}
\end{bmatrix},\qquad k\in\NN,
\qquad
\bX_{k,q}:=\begin{bmatrix}
X_{k,1}\\
\vdots\\
X_{k,q}
\end{bmatrix}, \qquad k\in\ZZ_+.
\end{split}
\end{gather}
Note that for each $i,j\in\{1,\dots,q\}$ and $r\in\{1,\dots,p\}$, the $(i,j)$-th entry of $\bV_q^{(r)}$ is
$v_{i,j}^{(r)}=\cov(\xi_{1,1,r,i},\xi_{1,1,r,j})$.
As it was explained at the beginning of Section \ref{Section_conv_results}, $\bA$ is a lower triangular matrix whose upper-left $q\times q$ submatrix is $\bI_q$,
thus, if $i\neq r$, $i\in\{1,\dots,q\}$, $r\in\{1,\dots,p\}$, then we have $a_{i,r}=0$, which implies $\PP(\xi_{1,1,r,i}=0)=1$ and hence $v_{i,j}^{(r)}=\cov(0,\xi_{1,1,r,j})=0$.
Replacing the role of $i$ by $j$ in the previous argument, if $j\neq r$, $j\in\{1,\dots,q\}$, $r\in\{1,\dots,p\}$, then we get that $v_{i,j}^{(r)}=0$.
Consequently, for each $r\in\{1,\dots,p\}$, all the entries of $\bV_q^{(r)}$ are equal to $0$, except possibly the $(r,r)$-th entry
 (which could happen only in case of $r\in\{1,\dots,q\}$).

Using \eqref{Mk1} together with the fact that $a_{i,j}=\delta_{i,j}$ for $i\in\{1,\dots,q\}$ and $j\in\{1,\dots,p\}$, we get that
\begin{equation*}
\bM_{k,q}=\bX_{k,q}-\bX_{k-1,q}-\bb_q, \qquad k\in\NN.
\end{equation*}
Consider the random step processes
 \begin{align}\label{help15}
   \bcM_t^{(n)} :=
\begin{bmatrix}
 \cM_{t,1}^{(n)} \\
\vdots \\
\cM_{t,q}^{(n)}  \end{bmatrix}
   := \frac{1}{n} \sum_{k=1}^\nt \bM_{k,q}
   = \frac{1}{n} \bX_{\nt,q} - \frac{\nt}{n} \bb_q , \qquad t \in \RR_+ , \qquad n \in \NN ,
 \end{align}
 where at the third equality we used that \ $\bX_0 = \bzero$.
We show that
 \begin{equation}\label{conv_bM}
  (\bcM_t^{(n)})_{t\in\RR_+} \distr (\bcM_t)_{t\in\RR_+} \qquad
  \text{as \ $n \to \infty$,}
 \end{equation}
 where the limit process $(\bcM_t)_{t\in\RR_+} = ([\cM_{t,1}, \dots, \cM_{t,q}]^\top)_{t\in\RR_+}$ is the pathwise unique strong solution of the SDE
 \begin{equation}\label{SDE_bM}
  \begin{cases}
   \dd \cM_{t,1}
   &= \sqrt{v^{(1)}_{1,1} (\cM_{t,1} + b_1 t)^+ } \, \dd \cW_{t,1},\\[1mm]
  & \vdots \\[1mm]
   \dd \cM_{t,q}
   &= \sqrt{v^{(q)}_{q,q} (\cM_{t,q} + b_q t)^+ } \, \dd \cW_{t,q},
  \end{cases}
  \qquad t \in\RR_+,
 \end{equation}
 with initial value $ \bcM_0 = \bzero\in\RR^q$, where $(\cW_{t,1})_{t\in\RR_+},\dots,(\cW_{t,q})_{t\in\RR_+}$
 are independent standard Wiener processes.

We check that the SDE \eqref{SDE_bM} admits a pathwise unique strong solution $(\bcM_{t})_{t\in\RR_+}$ for all initial values $\bx\in\RR^q$.
If $(\bcM_{t})_{t\in\RR_+}$ is a pathwise unique strong solution of the SDE \eqref{SDE_bM} with initial value $\bx\in\RR^q$,
 then, by It\^o's formula, $\bcP_t = [\cP_{t,1},\ldots,\cP_{t,q}]^\top :=\bcM_{t}+\bb_q t$, $t\in\RR_+$, is a pathwise unique strong solution of the SDE
 \begin{equation}\label{SDE_P}
  \begin{cases}
   &\dd\cP_{t,1}
   = b_1 \, \dd t
     + \sqrt{v^{(1)}_{1,1} \, \cP_{t,1}^+} \,
       \dd \cW_{t,1} , \\[1mm]
   &\phantom{\dd\cX_{t,1}\;\;}\vdots\\[1mm]
   &\dd\cP_{t,q}
   = b_{q} \, \dd t
     + \sqrt{v^{(q)}_{q,q} \, \cP_{t,q}^+ } \,
       \dd \cW_{t,q} ,
  \end{cases}
  \qquad t \in\RR_+ ,
 \end{equation}
with initial value $\bcP_0=\bx\in\RR^q$.
Conversely, if $(\bcP_t)_{t\in\RR_+}$ is a pathwise unique strong solution of the SDE \eqref{SDE_P} with initial value $\bcP_0=\bp\in\RR^q$,
then, by It\^o's formula, the process $\widetilde\bcM_{t} :=\bcP_t-\bb_q t$, $t\in\RR_+$,
 is a pathwise unique strong solution of the SDE \eqref{SDE_bM} with initial value $\bp\in\RR^q$.
Consequently, the SDE \eqref{SDE_bM} admits a pathwise unique strong solution if and only if
the SDE \eqref{SDE_P} admits such a solution as well.
 Using Ikeda and Watanabe \cite[Chapter IV, Example 8.2]{IkeWat},
 the SDE \eqref{SDE_P} admits a pathwise unique strong solution with all initial values $\bp\in\RR^q$.
 Further, since $\bcM_0=\bzero\in\RR^q$ and $b_i, v_{i,i}^{(i)}\in\RR_+$, $i\in\{1,\ldots,q\}$,
 we have that $\bcM_t+\bb_qt=\bcP_t\in\RR_+^q$ almost surely for all $t\in\RR_+$.

{\sl Step 2/(b).}
To prove \eqref{conv_bM}, we will use Theorem \ref{Conv2DiffThm} with the choices
$d=r=q$,
$\bcU=\bcM$,
$\bU_k^{(n)}= n^{-1}\bM_{k,q}$, $n,k\in\NN$,
$\bU_0^{(n)} =\bzero\in\RR^q$, $n\in\NN$,
$\cF_k^{(n)}=\cF_k^{\bX}$, $n\in\NN$, $k\in\ZZ_+$
(yielding $\bcU^{(n)}=\bcM^{(n)}$, $n\in\NN$),
and with coefficient functions $\bbeta:\RR_+\times\RR^q\to\RR^q$ and $\bgamma:\RR_+\times\RR^q\to\RR^{q\times q}$
given by $\bbeta(t,\bx):=\bzero$ and
\begin{equation*}
\bgamma(t,\bx):=
\begin{bmatrix}
\sqrt{v_{1,1}^{(1)}(x_1+b_1t)^+} & 0 & \dots & 0\\
0 & \sqrt{v_{2,2}^{(2)}(x_2+b_2t)^+} & \dots & 0\\
\vdots & \vdots & \ddots &\vdots\\
0 & 0 & \dots & \sqrt{v_{q,q}^{(q)}(x_q+b_qt)^+}
\end{bmatrix}
\end{equation*}
for all $t\in\RR_+$ and $\bx\in\RR^q$.
 With these notations, the SDE \eqref{SDE_bM} takes the form
\begin{equation*}
\dd\bcM_t=\bbeta(t,\bcM_t)\,\dd t+\bgamma(t,\bcM_t)\dd\bcW_t,\qquad t\in\RR_+,
\end{equation*}
with initial value $\bcM_0=\bzero\in\RR^q$.
The convergence $\bU_0^{(n)}\distr \bzero$ as $n\to\infty$, and
 (using that $\EE(\bM_{k,q} \mid \cF_{k-1}^\bX)=\bzero$, $k\in\NN$) condition (i) of Theorem \ref{Conv2DiffThm} trivially hold.
Condition (ii) of Theorem \ref{Conv2DiffThm} requires that we have
\begin{equation}\label{cond_ii}
\sup_{t\in[0,T]}
   \biggl\Vert\frac{1}{n^2}
          \sum_{k=1}^\nt
           \EE(\bM_{k,q} \bM_{k,q}^\top \mid \cF_{k-1}^\bX)
          - \int_0^t \bgamma(s,\bcM_s^{(n)})\bgamma(s,\bcM_s^{(n)})^\top \, \dd s \biggr\Vert
  \stoch 0 \qquad \text{as \ $n\to\infty$}
\end{equation}
for all $T\in\RR_+$, where, since $\bgamma(s,\bx)$, $s\in\RR_+$, $\bx\in\RR^q$, is diagonal, for all $s\in\RR_+$ and $n\in\NN$, we have
\begin{gather*}
  \bgamma(s,\bcM_s^{(n)})\bgamma(s,\bcM_s^{(n)})^\top
  = \begin{bmatrix}
      v^{(1)}_{1,1}(\cM_{s,1}^{(n)} + b_1 s)^+ & 0 & \dots & 0 \\
      0 & v^{(2)}_{2,2}(\cM_{s,2}^{(n)} + b_2 s)^+ & \dots & 0 \\
\vdots & \vdots & \ddots & \vdots\\
      0 & 0 & \dots &  v^{(q)}_{q,q}(\cM_{s,q}^{(n)} + b_q s)^+ \\
     \end{bmatrix}.
 \end{gather*}
Due to the form that the matrices $\bV_q^{(r)}$, $r\in\{1,\dots,q\}$ have (see Step 2/(a)), we can write this as
\begin{equation*}
\bgamma(s,\bcM_s^{(n)})\bgamma(s,\bcM_s^{(n)})^\top=\sum_{r=1}^q(\cM_{s,r}^{(n)}+b_r s)^+\bV_q^{(r)},\qquad s\in\RR_+,n\in\NN.
\end{equation*}
By \eqref{help15}, we have that
\begin{equation*}
   \bcM_s^{(n)} + \bb_q s
   = \frac{1}{n} \bX_{\ns,q} + \frac{ns-\ns}{n} \bb_q, \qquad n\in\NN, s\in\RR_+.
\end{equation*}
Since $(\bX_{k,q})_{k\in\ZZ_+}$ has non-negative coordinate processes and $\bb_q\in\RR_+^q$,
we get that $\bcM_s^{(n)} + \bb_q s\in\RR_+^q$, $s\in\RR_+$, $n\in\NN$.
Consequently, for $n\in\NN$, $t\in\RR_+$, we have
 \begin{align}\label{help_06}
 \begin{split}
 \int_0^t\bgamma(s,\bcM_s^{(n)})\bgamma(s,\bcM_s^{(n)})^\top\,\dd s
  &=\sum_{i=1}^q\left(\int_0^t(\cM_{s,i}^{(n)} + b_i s)^+\,\dd s\right)\bV_q^{(i)}\\
  &=\sum_{i=1}^q\left(\int_0^t(\cM_{s,i}^{(n)} + b_i s)\,\dd s\right)\bV_q^{(i)}.
 \end{split}
 \end{align}
 Note that, for each $k\in\NN$, $\bM_{k,q}\bM_{k,q}^\top$ is the upper-left $q\times q$ submatrix of $\bM_k\bM_k^\top$, and $\bV_q^{(r)}$ is the upper-left $q\times q$ submatrix of $\bV^{(r)}$, $r\in\{0,\dots,p\}$.
 Using \eqref{condVa2r} and the fact that $\bV_q^{(r)}=\bzero\in\RR^{q\times q}$, $r\in\{q+1,\dots,p\}$ (see Step 2/(a)),
we obtain that
 \begin{align*}
   \frac{1}{n^2}
   \sum_{k=1}^\nt
    \EE(\bM_{k,q} \bM_{k,q}^\top \mid \cF_{k-1}^\bX)
   &= \frac{\nt}{n^2} \bV_q^{(0)}
      + \frac{1}{n^2}
       \sum_{k=1}^\nt \left(\sum_{i=1}^p X_{k-1,i} \bV_q^{(i)}\right)\\
   & = \frac{\nt}{n^2} \bV_q^{(0)}
      + \frac{1}{n^2}
       \sum_{k=1}^\nt \left(\sum_{i=1}^q X_{k-1,i} \bV_q^{(i)}\right)
 \end{align*}
 for all $t \in \RR_+$ and $n \in \NN$.
This together with \eqref{help_06} implies that
\begin{align*}
&\frac{1}{n^2}
          \sum_{k=1}^\nt
           \EE(\bM_{k,q} \bM_{k,q}^\top \mid \cF_{k-1}^\bX)
          - \int_0^t \bgamma(s,\bcM_s^{(n)})\bgamma(s,\bcM_s^{(n)})^\top\, \dd s\\
&\qquad\qquad=\frac{\nt}{n^2} \bV_q^{(0)}
     + \sum_{i=1}^q\left(\frac{1}{n^2}
       \sum_{k=1}^\nt X_{k-1,i}-\int_0^t(\cM_{s,i}^{(n)} + b_i s)\,\dd s\right) \bV_q^{(i)}
\end{align*}
for all $n\in\NN$  and $t\in\RR_+$.
By the triangle inequality,
 \begin{align*}
  &\left\Vert \frac{\nt}{n^2} \bV_q^{(0)}
      + \sum_{i=1}^q\left(\frac{1}{n^2}
        \sum_{k=1}^\nt X_{k-1,i}-\int_0^t(\cM_{s,i}^{(n)} + b_i s)\,\dd s\right) \bV_q^{(i)}
     \right\Vert\\
  &\qquad \leq \frac{\nt}{n^2} \Vert \bV_q^{(0)} \Vert
                + \sum_{i=1}^q\left\vert \frac{1}{n^2} \sum_{k=1}^\nt X_{k-1,i}-\int_0^t(\cM_{s,i}^{(n)} + b_i s)\,\dd s \right\vert
                    \Vert \bV_q^{(i)} \Vert,
                \qquad n\in\NN,\; t\in\RR_+,
 \end{align*}
 and consequently, in order to show \eqref{cond_ii}, it is enough to verify that, for all $T\in\RR_{++}$, it holds that
 \[
  \sup_{t\in[0,T]} \frac{\nt}{n^2} \to 0 \qquad \text{as $n\to\infty$,}
 \]
 and
 \[
      \sup_{t\in[0,T]}  \left\vert \frac{1}{n^2} \sum_{k=1}^\nt X_{k-1,i}-\int_0^t(\cM_{s,i}^{(n)} + b_i s)\,\dd s \right\vert
        \stoch 0 \qquad \text{as $n\to\infty$,} \qquad i\in\{1,\ldots,q\}.
 \]
It is clear that
\begin{equation*}
\sup_{t\in[0,T]} \frac{\nt}{n^2}\leq \frac{T}{n}\to0 \qquad \text{as $n\to\infty$}
\end{equation*}
holds for all $T\in\RR_+$.
For the rest of this step, let $i\in\{1,\dots,q\}$ be fixed.
It is easy to see that (for the details, see, e.g., the proof of Theorem 1.1 in Barczy et al.\ \cite{BarBezPap1}), for $n\in\NN$ and $t\in\RR_+$, we have
\begin{equation*}
\int_0^t(\cM_{s,i}^{(n)} + b_i s)\,\dd s=\frac{1}{n^2}\sum_{k=0}^{\nt-1}X_{k,i}+\frac{nt-\nt}{n^2}X_{\nt,i}+\frac{\nt+(nt-\nt)^2}{2n^2}b_i,
\end{equation*}
 and thus
\begin{align*}
       \left|\frac{1}{n^2}\sum_{k=1}^\nt X_{k-1,i}-\int_0^t(\cM_{s,i}^{(n)} + b_i s)\,\dd s\right|
	&=\int_0^t(\cM_{s,i}^{(n)} + b_i s)\,\dd s -\frac{1}{n^2}\sum_{k=1}^\nt X_{k-1,i}\\
	&=\frac{nt-\nt}{n^2}X_{\nt,i}+\frac{\nt+(nt-\nt)^2}{2n^2}b_i\\
	&\leq n^{-2}X_{\nt,i}+\frac{\nt+1}{n^2}b_i.
\end{align*}
Since
\begin{equation*}
\sup_{t\in[0,T]} \frac{\nt+1}{n^2} \leq  \frac{T}{n}+\frac{1}{n^2} \to 0 \qquad \text{as $n\to\infty$}
\end{equation*}
for all $T\in\RR_+$, in order to finish the proof of \eqref{cond_ii}, it is enough to check that
\begin{equation*}
\sup_{t\in[0,T]} n^{-2}X_{\nt,i} \stoch0 \qquad \text{as $n\to\infty$}
\end{equation*}
for all $T\in\RR_+$.
To prove this, we verify that
\begin{equation*}
\EE\left(
  \frac{1}{n^4}
  \sup_{t \in [0,T]}
       X_{\nt,i}^2
\right)\to0 \qquad \text{as $n\to\infty$}
\end{equation*}
holds for all $T\in\RR_+$, which implies the desired convergence in probability.
Since $i\in\{1,\dots,q\}$, we have
 \begin{equation}\label{help_09}
    X_{k,i}=  X_{k,i}^{(0)} =\sum_{\ell=1}^k(M_{\ell,i}+b_i), \qquad k\in\NN,
 \end{equation}
 where the first equality follows by Step 1/(c) and the second equality by \eqref{M_sums}.
This, together with Lemma \ref{M_sum_growth} and $\eta_i=1$ (following from the definition of $q$) yields that
\begin{align*}
\EE\left(
  \frac{1}{n^4}
  \sup_{t \in [0,T]}
       X_{\nt,i}^2
\right)
 =
  \frac{1}{n^4}
\EE\left(
  \sup_{t \in [0,T]}
   \left(\sum_{\ell=1}^\nt(M_{\ell,i}+b_i)\right)^2
\right)
=\frac{1}{n^4}\OO(n^{1+1})=\OO(n^{-2})\to0
\end{align*}
as $n\to\infty$ for all $T\in\RR_+$.

Finally, we check condition (iii) of Theorem \ref{Conv2DiffThm}, which requires that
 \[
   \frac{1}{n^2}\sum_{k=1}^{\lfloor nT \rfloor}
     \EE\bigl(\|\bM_{k,q}\|^2 \bbone_{\{\|\bM_{k,q}\| > n\theta\}}
             \bmid \cF^{\bX}_{k-1}\bigr)
   \stoch 0 \qquad \text{as $n\to\infty$}
 \]
 for all $T\in\RR_{++}$ and $\theta\in\RR_{++}$.
Using that $\|\bM_{k,q}\|^2=\sum_{i=1}^q M_{k,i}^2$, $k\in\NN$, for each $n,k\in\NN$ and all $\theta\in\RR_{++}$,
 we have that
\[
\bbone_{\{\|\bM_{k,q}\|>n\theta\}}
	=\bbone_{\{\|\bM_{k,q}\|^2>n^2\theta^2\}}
	=\bbone_{\{\sum_{i=1}^q M_{k,i}^2>n^2\theta^2\}}
	\leq\sum_{i=1}^q\bbone_{\{M_{k,i}^2>n^2\theta^2/q\}}
	=\sum_{i=1}^q\bbone_{\{|M_{k,i}|>n\theta/\sqrt{q}\}}.
\]
Consequently, it is enough to show that
 \begin{equation}\label{alt_eq1}
   \frac{1}{n^2}\sum_{k=1}^{\lfloor nT \rfloor}
     \EE\bigl(M_{k,i}^2 \bbone_{\{|M_{k,i}| > n\theta\}}
             \bmid \cF^{\bX}_{k-1}\bigr)
   \stoch 0 \qquad \text{as $n\to\infty$, $i\in\{1,\dots,q\}$,}
 \end{equation}
and
 \begin{equation}\label{alt_eq2}
   \frac{1}{n^2}\sum_{k=1}^{\lfloor nT \rfloor}
     \EE\bigl(M_{k,i}^2 \bbone_{\{|M_{k,j}| > n\theta\}}
             \bmid \cF^{\bX}_{k-1}\bigr)
   \stoch 0 \qquad \text{as $n\to\infty$, $i\neq j$, $i,j\in\{1,\dots,q\}$}
 \end{equation}
 for all $T\in\RR_{++}$ and $\theta\in\RR_{++}$.
Here $M_{k,i}=X_{k,i}-X_{k-1,i}-b_i$, $k\in\NN$, $i\in\{1,\dots,q\}$ (following from \eqref{help_09}),
 where $(X_{k,i})_{k\in\ZZ_+}$ is a single-type critical GWI process 
 (following from the first paragraph of Section \ref{Section_conv_results} and the fact that the upper-left $q\times q$ submatrix of $\bA$ 
  is $\bI_q$, as it was explained at the beginning of Section \ref{Section_conv_results}).
Consequently, we can see that \eqref{alt_eq1} has already been shown in Step 3 of the proof 
 of Theorem 1.1 in Barczy et al.\,\cite{BarBezPap1}.
If $q=1$, then \eqref{alt_eq2} holds trivially.
For the rest of this part of the proof, suppose that $q>1$.
Using again that the upper-left $q\times q$ submatrix of $\bA$ is $\bI_q$, we can consider the decomposition
\[
M_{k,h}=N_{k,h}+(\vare_{k,h}-b_h), \qquad k\in\NN, \quad h\in\{1,\dots,q\},
\]
where
\[
N_{k,h}:=\sum_{\ell=1}^{X_{k-1,h}}(\xi_{k,\ell,h,h}-1),	\qquad k\in\NN, \quad h\in\{1,\dots,q\}.
\]
For the rest of this proof, let $T\in\RR_{++}$, $\theta\in\RR_{++}$, and $i\neq j$, $i,j\in\{1,\dots,q\}$ be fixed.
It is easy to see that for each $n,k\in\NN$, we have
\begin{equation}\label{alt_eq3}
M_{k,i}^2\bbone_{\{|M_{k,j}| > n\theta\}}
	\leq	2N_{k,i}^2\bbone_{\{|N_{k,j}| > n\theta/2\}}
		+ 2N_{k,i}^2\bbone_{\{|\vare_{k,j}-b_j| > n\theta/2\}}
		+ 2 (\vare_{k,i}-b_i)^2,
\end{equation}
for more details, see Barczy et al.\,\cite[Step 3 of the proof of Theorem 1.1]{BarBezPap1}.
Using that the random variables $\{\xi_{k,\ell,i,i}, \xi_{k,\ell,j,j}: \ell\in\NN\}$ are independent
 of the $\sigma$-algebra $\cF_{k-1}^\bX$, $k\in\NN$,
 the independence of $\{\xi_{k,\ell,i,i}: k,\ell\in\NN\}$ and $\{\xi_{k,\ell,j,j}: k,\ell\in\NN\}$,
 and Markov's inequality, we have that
\begin{align*}
	\EE\left(N_{k,i}^2\bbone_{\{|N_{k,j}| > n\theta\}}\mid\cF_{k-1}^\bX\right)
		&=\EE\left(\left(\sum_{\ell=1}^y(\xi_{k,\ell,i,i}-1)\right)^2\bbone_{\{ | \sum_{\ell=1}^{z}(\xi_{k,\ell,j,j}-1)| > n\theta\}}\right){\Bigg |}_{\,y=X_{k-1,i},\, z=X_{k-1,j}}\\
		&=X_{k-1,i}v_{i,i}^{(i)}\PP(|N_{k,j}| > n\theta\mid\cF_{k-1}^\bX)\\
		&\leq X_{k-1,i}v_{i,i}^{(i)}\EE(N_{k,j}^2\mid\cF_{k-1}^\bX)n^{-2}\theta^{-2}\\
		&=X_{k-1,i}v_{i,i}^{(i)}X_{k-1,j}v_{j,j}^{(j)}n^{-2}\theta^{-2}\\
		&\leq (X_{k-1,i}^2+X_{k-1,j}^2)v_{i,i}^{(i)}v_{j,j}^{(j)}n^{-2}\theta^{-2}, \qquad n,k\in\NN.
\end{align*}
Consequently,
\[
   \frac{1}{n^2}\sum_{k=1}^{\lfloor nT \rfloor}
	\EE\left(N_{k,i}^2\bbone_{\{|N_{k,j}| > n\theta\}}\mid\cF_{k-1}^\bX\right)
   \stoch 0 \qquad \text{as $n\to\infty$,}
\]
since
\[
n^{-4}\sum_{k=1}^\nT\EE(X_{k-1,i}^2+X_{k-1,j}^2)=n^{-4}\sum_{k=1}^\nT\OO(k^2)=n^{-4}\OO(n^3)=\OO(n^{-1}), \qquad n\in\NN,
\]
where we used that $\EE(X_{k,i}^2+X_{k,j}^2)=\OO(k^2)$, $k\in\ZZ_+$, by Remark \ref{Rem_X2}.
Furthermore, we can show that
\begin{equation}\label{alt_eq4}
  \frac{1}{n^2}\sum_{k=1}^{\lfloor nT \rfloor}
	\EE\left(N_{k,i}^2\bbone_{\{|\vare_{k,j}-b_j| > n\theta\}}\mid\cF_{k-1}^\bX\right)
   \stoch 0 \qquad \text{as $n\to\infty$,}
\end{equation}
 by a similar argument as in Barczy et al.\,\cite[Step 3b of the proof of Theorem 1.1]{BarBezPap1}.
Namely, using that the random variables $\{\xi_{k,\ell,i,i},\vare_{k,j}:\ell\in\NN\}$ are independent and independent of the $\sigma$-algebra $\cF_{k-1}^\bX$, $k\in\NN$, we have that
\begin{align*}
\EE\left(N_{k,i}^2\bbone_{\{|\vare_{k,j}-b_j| > n\theta\}}\mid\cF_{k-1}^\bX\right)
	&=\EE\left(\left(\sum_{\ell=1}^{z}(\xi_{k,\ell,i,i}-1)\right)^2\bbone_{\{|\vare_{k,j}-b_j| > n\theta\}}\right){\Bigg |}_{z=X_{k-1,i}}\\
	&=X_{k-1,i}v_{i,i}^{(i)}\PP(|\vare_{k,j}-b_j| > n\theta)\\
	&\leq	X_{k-1,i}v_{i,i}^{(i)}\EE((\vare_{k,j}-b_j)^2)n^{-2}\theta^{-2}\\
	&=X_{k-1,i}v_{i,i}^{(i)}v_{j,j}^{(0)}n^{-2}\theta^{-2}, \qquad n,k\in\NN.
\end{align*}
This implies that \eqref{alt_eq4} holds, since 
 \begin{align}\label{help_10}
  \frac{1}{n^4}\sum_{k=1}^\nT X_{k-1,i}\stoch 0 \qquad \text{as $n\to\infty$.}
 \end{align}
Indeed, $\EE(X_{k-1,i})=\OO(k)$, $k\in\NN$, by Lemma \ref{EEX1} and $\eta_i=1$ (since $i\leq q$),
 and hence 
 \[
  \EE\left( \frac{1}{n^4}\sum_{k=1}^\nT X_{k-1,i} \right)
     = \frac{1}{n^4} \sum_{k=1}^\nT \OO(k) 
     = \OO(n^{-2}),\qquad n\in\NN,
 \]
 yielding \eqref{help_10}.
Finally, the convergence
\[
 \frac{1}{n^2}\sum_{k=1}^{\lfloor nT \rfloor}
	\EE\left((\vare_{k,i}-b_i)^2\mid\cF_{k-1}^\bX\right)
   \stoch 0 \qquad \text{as $n\to\infty$}
\]
 has already been shown in Barczy et al.\,\cite[Step 3c of the proof of Theorem 1.1]{BarBezPap1}.
Altogether, the previous considerations imply that condition (iii) of Theorem \ref{Conv2DiffThm} is satisfied,
 and hence we finished the proof of \eqref{conv_bM}.

{\sl Step 2/(c).} Now we show that $(n^{-1}\bX_{\nt,q})_{t\in\RR_+}$ converges weakly
  towards $([\cX_{t,1},\ldots,\cX_{t,q}]^\top)_{t\in\RR_+}$ as $n\to\infty$, where
 $\cX_{t,i}:=\cM_{t,i} + tb_i$, $t\in\RR_+$, $i\in\{1,\ldots,q\}$.
Note that, by \eqref{help15}, using the notation for the maps $\tau_n$, $n\in\NN$, defined in part (iv) of Definition \ref{maps}, we have
 \begin{equation}\label{help_02}
 \frac{1}{n}\bX_{\nt,q}=\bcM_t^{(n)}+\frac{\nt}{n}\bb_q
 =\left(\tau_n\left(
 \left( \bcM_s^{(n)} + s\bb_q \right)_{s\in\RR_+}
 \right)\right)(t), \qquad n\in\NN,\quad t\in\RR_+,
 \end{equation}
 where we used the fact that $\left(\tau_n\left((\bcM_s^{(n)})_{s\in\RR_+}\right)\right)(t)=\bcM_t^{(n)}$,
 $n\in\NN$, $t\in\RR_+$ (due to $\frac{\lfloor n \frac{\nt}{n}\rfloor}{n}=\frac{\nt}{n}$).
 By \eqref{conv_bM}, we have that $(\bcM_t^{(n)})_{t\in\RR_+} \distr (\bcM_t)_{t\in\RR_+}$ as $n\to\infty$,
 and clearly, the function $(t\bb_q)_{t\in\RR_+}$ belongs to $\CC(\RR_+,\RR^q)$.
Then, by Jacod and Shiryaev \cite[Chapter VI, Proposition 3.17]{JacShi},
 we get that
\begin{equation*}
 \left(\bcM_t^{(n)}+t\bb_q \right)_{t\in\RR_+}
  \distr
 \left(\bcM_t + t\bb_q \right)_{t\in\RR_+}  \qquad  \text{as $n\to\infty$.}
\end{equation*}
Since $(\bcM_t)_{t\in\RR_+}$ has continuous trajectories almost surely,
  using \eqref{help_02} and Lemma \ref{time_transform}, we obtain that
 \begin{equation*}
 \left(n^{-1}\bX_{\nt,q}\right)_{t\in\RR_+}=\left(\bcM_t^{(n)}+\frac{\nt}{n}\bb_q\right)_{t\in\RR_+}
  \distr(\bcM_t+t\bb_q)_{t\in\RR_+}
           = \left(\begin{bmatrix}
                    \cX_{t,1} \\
                    \vdots \\
                    \cX_{t,q} \\
                  \end{bmatrix}\right)_{t\in\RR_+}
          \qquad \text{as $n\to\infty$,}
 \end{equation*}
 where the limit process $([\cX_{t,1},\ldots,\cX_{t,q}]^\top)_{t\in\RR_+}$
  is the pathwise unique strong solution of the SDE \eqref{SDE_X_1} (which can be checked by It\^o's formula).

\smallskip

\textsl{Step 3.}  Asymptotic behavior of $(\bX_k)_{k\in\ZZ_+}$.

{\sl Step 3/(a).} Now, we turn to prove \eqref{Conv_X_1}.
Using \eqref{decomp_vector}, for $n\in\NN$, we can write
\begin{align}\label{decomp_vec2}
\begin{split}
&\left(\begin{bmatrix}
n^{-\eta_1}X_{\nt,1}\\
\vdots\\
n^{-\eta_p}X_{\nt,p}
\end{bmatrix}\right)_{t\in\RR_+}\\
&=\left(\begin{bmatrix}
n^{-\eta_1} \sum_{m=0}^{\eta_1-2}\sum_{j=1}^{p}c_{1,j}^{[m]}X_{\nt,j}^{(m)}\\
\vdots\\
n^{-\eta_p} \sum_{m=0}^{\eta_p-2}\sum_{j=1}^{p}c_{p,j}^{[m]}X_{\nt,j}^{(m)}
\end{bmatrix}\right)_{t\in\RR_+}
+
\left(\begin{bmatrix}
n^{-\eta_1}\sum_{j=1}^{q}c_{1,j}^{[\eta_1-1]}X_{\nt,j}^{(\eta_1-1)}\\
\vdots\\
n^{-\eta_p}\sum_{j=1}^{q}c_{p,j}^{[\eta_p-1]}X_{\nt,j}^{(\eta_p-1)}
\end{bmatrix}\right)_{t\in\RR_+}.
\end{split}
\end{align}
We will proceed by handling the two stochastic processes at the right-hand side of \eqref{decomp_vec2} separately, starting with the second one.
Namely, we show that the second stochastic process at the right-hand side of \eqref{decomp_vec2}
converges weakly to $\left(\bcX_t\right)_{t\in\RR_+}$ as $n\to\infty$, where $\left(\bcX_t\right)_{t\in\RR_+}$ is the limit process in \eqref{Conv_X_1}.
Note that, by \eqref{integral_dec}, for $i,j\in\{1,\dots,p\}$,
 due to the linearity of integration, we have
 \begin{align*}
 n^{-\eta_i}X_{\nt,j}^{(\eta_i-1)}
	=&n^{-\eta_i}n^{\eta_i-1}\left(\varphi_n^{\circ(\eta_i-1)}\left(\left(X_{\nr,j}^{(0)}\right)_{r\in\RR_+}\right)\right)(t)\\
	=&\left(\varphi_n^{\circ(\eta_i-1)}\left(\left(n^{-1}X_{\nr,j}^{(0)}\right)_{r\in\RR_+}\right)\right)(t), \qquad n\in\NN, t\in\RR_+,
\end{align*}
 where we used the notation $\varphi_n^{\circ(m)}$, $n,m\in\NN$, defined in Definition \ref{maps}.
 Further, by \eqref{d_int}, for $i\in\{1,\dots,p\}$ and $j\in\{1,\dots,q\}$, we have
\begin{equation*}
 n^{-\eta_i}X_{\nt,j}^{(\eta_i-1)}=\left(\varphi_n^{\circ(\eta_i-1)}\left(\left(n^{-1}X_{\nr,j}\right)_{r\in\RR_+}\right)\right)(t), \qquad n\in\NN, t\in\RR_+.
\end{equation*}

Consequently, we have
\begin{align}\label{help_01}
 \begin{bmatrix}
   n^{-\eta_1}\sum_{j=1}^qc_{1,j}^{[\eta_1-1]}X_{\nt,j}^{(\eta_1-1)}\\
   \vdots\\
    n^{-\eta_p}\sum_{j=1}^qc_{p,j}^{[\eta_p-1]}X_{\nt,j}^{(\eta_p-1)}
  \end{bmatrix}
 =
 \begin{bmatrix}
  \sum_{j=1}^qc_{1,j}^{[\eta_1-1]} \left(\varphi_{n}^{\circ(\eta_1-1)}\left(\left(n^{-1}X_{\nr,j}\right)_{r\in\RR_+}\right)\right)(t)\\
  \vdots\\
  \sum_{j=1}^qc_{p,j}^{[\eta_p-1]}\left(\varphi_{n}^{\circ(\eta_p-1)}\left(\left(n^{-1}X_{\nr,j}\right)_{r\in\RR_+}\right)\right)(t)
 \end{bmatrix}
\end{align}
 for $n\in\NN$ and $t\in\RR_+$.

We check that the right-hand side of \eqref{help_01} can be written as
\begin{equation*}
\left(\bB\Phi_n^{(p-1)}\left(\left(
\begin{bmatrix}
n^{-1}X_{\nr,1}\\
\vdots\\
n^{-1}X_{\nr,q}
\end{bmatrix}
\right)_{r\in\RR_+}\right)\right)(t), \qquad n\in\NN, t\in\RR_+,
\end{equation*}
where $\Phi_n^{(p-1)}$ is defined in part (vii) of Definition \ref{maps} and the entries $b_{i,j}$, $i\in\{1,\ldots,p\}$, $j\in\{1,\ldots,qp\}$ of
 $\bB\in\RR^{p\times qp}$ are given by
\begin{align*}
b_{i,j}:=&\begin{cases}
c_{i,r}^{[\eta_i-1]} & \text{if $j=(\eta_i-1)q+r$ for some $r\in\{1,\dots,q\}$,}\\
0 & \text{otherwise}
\end{cases}\\
	=&\begin{cases}
c_{i,j-(\eta_i-1)q}^{[\eta_i-1]} & \text{if $j\in\{(\eta_i-1)q+1, \dots, \eta_iq\}$,}\\
0 & \text{otherwise},
\end{cases}
\end{align*}
 that is, for $i\in\{1,\dots,p\}$, the $i$-th row of $\bB$ is
  $\begin{bmatrix}
   0 & \dots & 0 & c_{i,1}^{[\eta_i-1]} & \dots & c_{i,q}^{[\eta_i-1]} &  0 & \dots & 0
  \end{bmatrix}$,
where we have $(\eta_i-1)q$ zeroes in the beginning and $(p-\eta_i)q$ zeroes at the end.
Indeed, for $n\in\NN$, $t\in\RR_+$, and $i\in\{1,\dots,p\}$, we have
\begin{align}\label{help_03}
\begin{split}
 \begin{bmatrix}b_{i,1} & \dots & b_{i, qp}\end{bmatrix}
 &\left(\Phi_n^{(p-1)}\left(\left(
 \begin{bmatrix}
   n^{-1}X_{\nr,1}\\
   \vdots\\
   n^{-1}X_{\nr,q}
 \end{bmatrix}
  \right)_{r\in\RR_+}\right)\right)(t)\\
 &=\begin{bmatrix}c_{i,1}^{[\eta_i-1]} & \dots & c_{i,q}^{[\eta_i-1]}\end{bmatrix}
  \left(\varphi_n^{\circ(\eta_i-1)}\left(\left(
   \begin{bmatrix}
   n^{-1}X_{\nr,1}\\
  \vdots\\
   n^{-1}X_{\nr,q}
 \end{bmatrix}
 \right)_{r\in\RR_+}\right)\right)(t)\\
 & =\begin{bmatrix}c_{i,1}^{[\eta_i-1]} & \dots & c_{i,q}^{[\eta_i-1]}\end{bmatrix}
          \begin{bmatrix}
             \left(  \varphi_n^{\circ(\eta_i-1)}\big( (n^{-1} X_{\nr,1})_{r\in\RR_+} \big) \right)(t) \\
            \vdots \\
             \left(  \varphi_n^{\circ(\eta_i-1)}\big( (n^{-1} X_{\nr,q})_{r\in\RR_+} \big) \right)(t) \\
          \end{bmatrix} \\
&=\sum_{j=1}^qc_{i,j}^{[\eta_i-1]}\left(\varphi_{n}^{\circ(\eta_i-1)}\left(\left(n^{-1}X_{\nr,j}\right)_{r\in\RR_+}\right)\right)(t),
\end{split}
\end{align}
as desired.
By Step 2 of this proof, we have that
\begin{equation*}
\left(
\begin{bmatrix}
n^{-1}X_{\nt,1}\\
\vdots\\
n^{-1}X_{\nt,q}
\end{bmatrix}
\right)_{t\in\RR_+}
\distr
\left(
\begin{bmatrix}
\cX_{t,1}\\
\vdots\\
\cX_{t,q}
\end{bmatrix}
\right)_{t\in\RR_+}
\qquad \text{as $n\to\infty$,}
\end{equation*}
where the limit process is the pathwise unique strong solution of the SDE \eqref{SDE_X_1}.
Then, as this limit process has continuous sample paths almost surely,
 by applying Lemmas \ref{ConvMatrixMult} and \ref{ConvIteratedIntegral},
we have that
\begin{equation*}
\left(\left(\bB\Phi_n^{(p-1)}\left(\left(
\begin{bmatrix}
n^{-1}X_{\nr,1}\\
\vdots\\
n^{-1}X_{\nr,q}
\end{bmatrix}
\right)_{r\in\RR_+}\right)\right)(t)\right)_{t\in\RR_+}
\distr
\left(\left(\bB\Phi^{(p-1)}\left(\left(
\begin{bmatrix}
\cX_{r,1}\\
\vdots\\
\cX_{r,q}
\end{bmatrix}
\right)_{r\in\RR_+}\right)\right)(t)\right)_{t\in\RR_+}
\end{equation*}
as $n\to\infty$.
Consequently, using \eqref{help_01} and the definitions of the matrix $\bB$ and
 the maps $\Phi^{(p-1)}$, $\Phi_{n}^{(p-1)}$, $n\in\NN$ (given in part (vii) in Definition \ref{maps}),
 by a similar calculation as in \eqref{help_03}, we have
\begin{equation}\label{help_04}
\left(\begin{bmatrix}
n^{-\eta_1}\sum_{j=1}^qc_{1,j}^{[\eta_1-1]}X_{\nt,j}^{(\eta_1-1)}\\
\vdots\\
n^{-\eta_p}\sum_{j=1}^qc_{p,j}^{[\eta_p-1]}X_{\nt,j}^{(\eta_p-1)}
\end{bmatrix}\right)_{t\in\RR_+}
\!\!\!\!\distr
\left(\begin{bmatrix}
\sum_{j=1}^qc_{1,j}^{[\eta_1-1]}\left(\varphi^{\circ(\eta_1-1)}\left(\left(\cX_{r,j}\right)_{r\in\RR_+}\right) \right)(t)\\
\vdots\\
\sum_{j=1}^qc_{p,j}^{[\eta_p-1]}\left(\varphi^{\circ(\eta_p-1)}\left(\left(\cX_{r,j}\right)_{r\in\RR_+}\right)\right)(t)
\end{bmatrix}
\right)_{t\in\RR_+}
\end{equation}
as $n \to \infty$.
The limit process in \eqref{help_04} coincides with $([\cX_{t,1},\ldots,\cX_{t,p}]^\top)_{t\in\RR_+}$
introduced in \eqref{SDE_X_1} and \eqref{SDE_X_2}.
Indeed, if $i\in\{1,\dots,q\}$, then we have $\eta_i=1$, and hence
 $c_{i,j}^{[\eta_i-1]} = c_{i,j}^{[0]} =\delta_{i,j}$ for $i,j\in\{1,\ldots,q\}$.
This together with the fact that $\varphi^{\circ(\eta_i-1)} = \varphi^{\circ(0)}$ is the identity mapping for $i\in\{1,\ldots,q\}$,
 implies that
 $\sum_{j=1}^q c_{i,j}^{[\eta_i-1]}\varphi^{\circ(\eta_i-1)}\big(\big(\cX_{r,j}\big)_{r\in\RR_+}\big)(t)
   = \cX_{t,i}$ for $i\in\{1,\ldots,q\}$ and $t\in\RR_+$.
Further, using the definition of $\varphi$ (see part (v) of Definition \ref{maps}),
 we readily have the first equality in \eqref{SDE_X_2}, that is, for each $i\in\{q+1,\dots,p\}$,
 the $i$-th coordinate process of the limit process at the right-hand side of \eqref{help_04} coincides
 with $(\cX_{t,i})_{t\in\RR_+}$.

Next, we show that for each $i\in\{q+1,\dots,p\}$ and $m\in\{1,\dots,\eta_i-1\}$, we have
\begin{equation}\label{intsum_change}
\sum_{j=1}^qc_{i,j}^{[\eta_i-1]}\left(\varphi^{\circ(\eta_i-1)}\left(\left(\cX_{r,j}\right)_{r\in\RR_+}\right)\right)(t)
 = \sum_{ \substack{j\in\{1,\dots,p\},\\\eta_j=\eta_i-m}} c_{i,j}^{[m]}\left(\varphi^{\circ(m)}\left(\left(\cX_{r,j}\right)_{r\in\RR_+}\right)\right)(t), \qquad t\in\RR_+.
\end{equation}
Let $i\in\{q+1,\dots,p\}$ and $m\in\{1,\dots,\eta_i-1\}$ be fixed.
Recall that for $r\in\{1,\dots,p\}$, if $c_{i,r}^{[m]}>0$,
then, by Lemma \ref{eta_ineq}, we have $\eta_i\geq\eta_r+m$.
On the other hand, by the definition of $\eta_r$, if $\eta_i\geq\eta_r+m+1$, then we have that
 $c_{r,j}^{[\eta_i-m-1]}=0$ for each $j\in\{1,\dots,p\}$.
Consequently, if $c_{i,r}^{[m]}c_{r,j}^{[\eta_i-m-1]}>0$ for some $r,j\in\{1,\ldots,p\}$, then $\eta_r=\eta_i-m$.
Hence we get that
\begin{align*}
\sum_{j=1}^qc_{i,j}^{[\eta_i-1]}\left(\varphi^{\circ(\eta_i-1)}\left((\cX_{s,j})_{s\in\RR_+}\right)\right)(t)
	=&\sum_{j=1}^q\sum_{r=1}^pc_{i,r}^{[m]}c_{r,j}^{[\eta_i-m-1]}\left((\varphi^{\circ(m)}\circ\varphi^{\circ(\eta_i-m-1)})\left((\cX_{s,j})_{s\in\RR_+}\right)\right)(t)\\
	=&\sum_{r=1}^p c_{i,r}^{[m]}  \left(\varphi^{\circ(m)} \left(\sum_{j=1}^q  c_{r,j}^{[\eta_i-m-1]} \varphi^{\circ(\eta_i-m-1)}\left((\cX_{s,j})_{s\in\RR_+}\right)\right)\right)(t)\\
	=&\sum_{r=1}^p  c_{i,r}^{[m]} \delta_{\eta_r,\eta_i-m}  \left(  \varphi^{\circ(m)} \left(\sum_{j=1}^q c_{r,j}^{[\eta_r-1]}  \varphi^{\circ(\eta_r-1)} \left((\cX_{s,j})_{s\in\RR_+}\right)\right)\right)(t)
\end{align*}
\begin{align*}
	\hspace{4.7cm}=&\sum_{r=1}^p  c_{i,r}^{[m]} \delta_{\eta_r,\eta_i-m}  \left( \varphi^{\circ(m)}  \left(\left(\cX_{s,r}\right)_{s\in\RR_+}\right)\right)(t)\\
	\hspace{4.7cm}=&\sum_{\substack{r\in\{1,\dots,p\}\\\eta_r=\eta_i-m}} c_{i,r}^{[m]} \left( \varphi^{\circ(m)} \left(\left(\cX_{s,r}\right)_{s\in\RR_+}\right)\right)(t),\qquad t\in\RR_+,
\end{align*}
 where the last but one step follows by the form of the limit process $([\cX_{t,1},\ldots,\cX_{t,p}]^\top)_{t\in\RR_+}$
 given in \eqref{help_04}.
This implies \eqref{intsum_change}, and consequently, we have the second equality in \eqref{SDE_X_2}.

{\sl Step 3/(b).} Next, we show that the first stochastic process at the right-hand side of \eqref{decomp_vec2} disappears in the limit.
More precisely, we show that for all $T\in\RR_+$,
\begin{equation}\label{stoch0}
\sup_{t\in[0,T]}\left\|
\begin{bmatrix}
n^{-\eta_1}\sum_{m=0}^{\eta_1-2}\sum_{j=1}^{p}c_{1,j}^{[m]}X_{\nt,j}^{(m)}\\
\vdots\\
n^{-\eta_p}\sum_{m=0}^{\eta_p-2}\sum_{j=1}^{p}c_{p,j}^{[m]}X_{\nt,j}^{(m)}
\end{bmatrix}
\right\|\stoch 0 \qquad \text{as \ $n \to \infty$.}
\end{equation}
Note that for all $T\in\RR_+$ and $n\in\NN$, we have
 \begin{align*}
   \sup_{t\in[0,T]}\left\|
   \begin{bmatrix}
     n^{-\eta_1}\sum_{m=0}^{\eta_1-2}\sum_{j=1}^{p}c_{1,j}^{[m]}X_{\nt,j}^{(m)}\\
      \vdots\\
     n^{-\eta_p}\sum_{m=0}^{\eta_p-2}\sum_{j=1}^{p}c_{p,j}^{[m]}X_{\nt,j}^{(m)}
    \end{bmatrix}
  \right\|
  \leq \sum_{i=1}^p \sup_{t\in[0,T]}
     \left\vert n^{-\eta_i}\sum_{m=0}^{\eta_i-2}\sum_{j=1}^{p}c_{i,j}^{[m]}X_{\nt,j}^{(m)}  \right\vert.
 \end{align*}
Hence, by the triangle inequality, it is enough to show that
\begin{equation}\label{help_05}
\sup_{t\in[0,T]}\left|
n^{-\eta_i}c_{i,j}^{[m]}X_{\nt,j}^{(m)}
\right|\stoch 0\qquad  \text{as $n\to\infty$}
\end{equation}
 for all $T\in\RR_+$ and for each $i,j\in\{1,\dots,p\}$ and $m\in\{0,\dots,\eta_i-2\}$.
We will prove that the convergence holds even in the $L_2$-norm.

In what follows, let $T\in\RR_+$, $i,j\in\{1,\dots,p\}$ and $m\in\{0,\dots,\eta_i-2\}$ be fixed.
Clearly, if $c_{i,j}^{[m]}=0$, then \eqref{help_05} holds.
By Lemma \ref{M_sum_growth} and \eqref{M_sums}, we have that
\begin{equation*}
\EE\left(\sup_{t\in[0,T]}\left(X_{\nt,j}^{(m)}\right)^2\right)=\EE\left(\sup_{t\in[0,T]}\left(\sum_{\ell=1}^\nt\binom{\nt-\ell}{m}(M_{\ell,j}+b_j)\right)^2\right)=\OO(n^{\eta_j+1+2m}), \qquad n\in\NN.
\end{equation*}
Consequently,
\begin{equation*}
\EE\left(\left(\sup_{t\in[0,T]}\left|n^{-\eta_i}X_{\nt,j}^{(m)}\right|\right)^2\right)=n^{-2\eta_i}\EE\left(\sup_{t\in[0,T]}\left(X_{\nt,j}^{(m)}\right)^2\right)=\OO(n^{\eta_j+1+2m-2\eta_i}), \qquad n\in\NN.
\end{equation*}
Recall that, by Lemma \ref{eta_ineq}, if $c_{i,j}^{[m]}>0$, then $\eta_i\geq\eta_j+m$,
thus whenever $c_{i,j}^{[m]}>0$ and $m\leq\eta_i-2$, we have
\begin{align*}
\eta_j+1+2m-2\eta_i
		=&(\eta_j+m-\eta_i)+(1+m-\eta_i)\leq0+(1+m-\eta_i)\\
		=&1+m-\eta_i\leq1+(\eta_i-2)-\eta_i=-1,
\end{align*}
 yielding that
\begin{equation*}
\EE\left(\left(\sup_{t\in[0,T]}\left|n^{-\eta_i}X_{\nt,j}^{(m)}\right|\right)^2\right)=\OO(n^{-1})\to0 \qquad \text{as \ $n \to \infty$,}
\end{equation*}
whenever $c_{i,j}^{[m]}>0$ and $m\leq\eta_i-2$.
This implies \eqref{help_05} whenever $c_{i,j}^{[m]}>0$ and $m\leq\eta_i-2$ as well,
thus we obtain that \eqref{stoch0} holds for all $T\in\RR_+$.

By Lemma \ref{Lem_JacShi} (a kind of Slutsky's lemma for stochastic processes with trajectories in $\DD(\RR_+, \RR^d)$),
 \eqref{decomp_vec2}, \eqref{help_04}, and \eqref{stoch0} imply \eqref{Conv_X_1}.
\proofend

\vspace*{5mm}

\appendix

\vspace*{5mm}

\noindent{\bf\Large Appendices}

\section{Some special sums}\label{appInt}

We present some auxiliary results for sums of (binomially weighted) values of a function defined on $\ZZ_+$ or $\NN$.
These results are used throughout the proofs.

\begin{Lem}\label{specSums}
Let $f:\ZZ_+\to\RR$ and $g:\NN\to\RR$.
Then, for each $k\in\ZZ_+$, $m,n\in\NN$, we have
\begin{align}\label{specSums_formula1}
 &\sum_{\ell=0}^k f(\ell) = n \int_{0}^{ \frac{k+1}{n} } f(\ns)\,\dd s,\\ \label{specSums_formula2}
 &\sum_{\ell=1}^k\binom{k-\ell}{m}g(\ell)=n\int_{0}^{\frac{k}{n}}\sum_{\ell=1}^\ns\binom{\ns-\ell}{m-1}g(\ell)\,\dd s,
\end{align}
 where recall the convention $\binom{p}{q}:=0$ for $p<q$, $p,q\in\ZZ_+$.
Consequently, for each $n,m\in\NN$ and $t\in\RR_+$, we get
\begin{align}\label{sumInt}
\begin{split}
&\sum_{\ell=1}^\nt\binom{\nt-\ell}{m}g(\ell)\\
	&=n^m\int_{0}^{\frac{\nt}{n}}\left(\int_0^{\frac{\lfloor n s_m\rfloor}{n}}
		\dots\left(\int_0^{\frac{\lfloor n s_3\rfloor}{n}}\left(
		\int_0^{\frac{\lfloor n s_2\rfloor}{n}}\sum_{\ell=1}^{\lfloor n s_1\rfloor}g(\ell)\,\dd s_1
		\right)\,\dd s_2\right)
		\dots\right)\,\dd s_m.
\end{split}
\end{align}
\end{Lem}

\noindent\textbf{Proof.}
The statement \eqref{specSums_formula1} is trivial for $k=0$, since then both sides of the equality are equal to $f(0)$.
The statement \eqref{specSums_formula2} is also trivial for $k=0$, since, due to the convention $\sum_{j=1}^0=0$,
 both sides of the equality are equal to $0$.
To prove \eqref{specSums_formula1}, let $n,k\in\NN$ be fixed.
Then we have that
 \begin{align*}
 \sum_{\ell=0}^k f(\ell)
    =n \sum_{\ell=0}^k f(\ell)\cdot \frac{1}{n}
    =n \sum_{\ell=0}^k f(\ell) \int_{\frac{\ell}{n}}^{\frac{\ell+1}{n}} 1\,\dd s
   =n \sum_{\ell=0}^k \int_{\frac{\ell}{n}}^{\frac{\ell+1}{n}} f(\ns)\,\dd s
   = n\int_{0}^{\frac{k+1}{n}} f(\ns)\,\dd s,
 \end{align*}
 yielding \eqref{specSums_formula1}.

To prove \eqref{specSums_formula2}, let $k,m,n\in\NN$ be fixed.
By \eqref{specSums_formula1}, we have that
\begin{align*}
  n\int_{0}^{\frac{k}{n}}\sum_{\ell=1}^\ns\binom{\ns-\ell}{m-1}g(\ell)\,\dd s
  	  =\sum_{h=0}^{k-1}\sum_{\ell=1}^h\binom{h-\ell}{m-1}g(\ell).
\end{align*}
Indeed, if $k\geq 2$, $k\in\NN$, then \eqref{specSums_formula1} readily implies the above equality,
 and if $k=1$, then it holds trivially, since then both of its two sides are equal to $0$ due to
 the convention $\sum_{j=1}^0:=0$.
Using again the convention $\sum_{j=1}^0:=0$ and changing the order of summations, we get that
\begin{align*}
&\sum_{h=0}^{k-1}\sum_{\ell=1}^h\binom{h-\ell}{m-1}g(\ell)
	    =\sum_{h=1}^{k-1}\sum_{\ell=1}^h\binom{h-\ell}{m-1}g(\ell)
	    =\sum_{\ell=1}^{k-1}\sum_{h=\ell}^{k-1}\binom{h-\ell}{m-1}g(\ell)\\
&\qquad =\sum_{\ell=1}^{k-1}\sum_{h=0}^{k-1-\ell}\binom{h}{m-1}g(\ell)
        =\sum_{\ell=1}^{k-1}\left(\sum_{h=m-1}^{k-\ell-1} \binom{h}{m-1}\right)g(\ell)
        =\sum_{\ell=1}^{k}\binom{k-\ell}{m}g(\ell),
\end{align*}
where the last but one equality holds due to the convention $\binom{p}{q}:=0$ for $p<q$, $p,q\in\ZZ_+$,
 and the last equality holds due to the well-known hockey-stick identity
\begin{equation*}
  \sum_{j=r}^{u}\binom{j}{r}=\binom{u+1}{r+1}, \qquad u,r\in\ZZ_+.
\end{equation*}
This yields \eqref{specSums_formula2}.

To prove \eqref{sumInt}, let $n\in\NN$ and $t\in\RR_+$  be fixed.
In the case of $m=1$, \eqref{sumInt} is equivalent to \eqref{specSums_formula2}, and therefore it holds.
Now let $m_0\in\NN$ be fixed, and suppose that \eqref{sumInt} holds when $m$ is replaced by $m_0$.
By \eqref{specSums_formula2}, we have
\begin{equation*}
\sum_{\ell=1}^\nt\binom{\nt-\ell}{m_0+1}g(\ell)
	=n\int_{0}^{\frac{\nt}{n}}\left(\sum_{\ell=1}^\ns\binom{\ns-\ell}{m_0}g(\ell)\right)\,\dd s,
\end{equation*}
which, together with the induction hypothesis, implies that
\begin{align*}
&\sum_{\ell=1}^\nt\binom{\nt-\ell}{ m_0+1}g(\ell)
	=n\int_{0}^{\frac{\nt}{n}}\left(\sum_{\ell=1}^\ns\binom{\ns-\ell}{m_0}g(\ell)\right)\,\dd s\\
&\qquad\qquad	
	=n\int_{0}^{\frac{\nt}{n}}
	\left(
	n^{m_0}\int_0^{\frac{\ns}{n}}
	\left(\int_0^{\frac{\lfloor n s_{m_0}\rfloor}{n}}
	\dots\left(\int_0^{\frac{\lfloor n s_2\rfloor}{n}}\sum_{\ell=1}^{\lfloor n s_1\rfloor}g(\ell)\,\dd s_1
	\right)\dots
	\right)\,\dd s_{m_0}\right)\,\dd s\\
&\qquad\qquad
	=n^{m_0+1}\int_{0}^{\frac{\nt}{n}}\left(
		\int_0^{\frac{\lfloor ns_{m_0+1}\rfloor}{n}}
		\left(
		\int_0^{\frac{\lfloor ns_{m_0}\rfloor}{n}}
		\dots\left(
		\int_0^{\frac{\lfloor n s_2\rfloor}{n}}\sum_{\ell=1}^{\lfloor n s_1\rfloor}g(\ell)\,\dd s_1
		\right)
		\dots\right)\,\dd s_{m_0}\right)\,\dd s_{m_0+1},
\end{align*}
 as desired.
\proofend

Next, we recall a result about the powers of certain lower triangular matrices, which plays a crucial role in the proofs.
The next lemma is only proven in the case $k\in\NN$ in Barczy and Bezd\'any \cite[Lemma 4.1]{BarBez},
 but it is easy to see that it holds for $k=0$ as well, since $\binom{0}{m}=0$ for each $m\in\{1,\dots,p-1\}$.
\begin{Lem}\label{matrixlemma}
Let $p\in\NN$ and let $\bA\in\RR^{p\times p}$
 be a lower triangular matrix such that $a_{i,i}=1$, $i\in\{1,\ldots,p\}$.
Then for each $k\in\ZZ_+$, we have
\begin{equation}\label{mPow}
\bA^k=\sum_{m=0}^{p-1}\binom{k}{m}\left(\bA-\bI_p\right)^m.
\end{equation}
\end{Lem}
In the next remark, we collect a few facts about the powers of $\bC:=\bA-\bI_p$ (see Barczy and Bezd\'any \cite[Remark 4.1]{BarBez}).
\begin{Rem}\label{matrix_power_remark}
Let $p\in\NN$ and let $\bA\in\RR_+^{p\times p}$ be a lower triangular matrix
such that $a_{i,i}=1$, $i\in\{1,\dots,p\}$, and recall that $\bC=\bA-\bI_p$.
Then, for each $i,j\in\{1,\dots,p\}$, we have
\begin{equation}\label{C_0}
  c_{i,j}^{[m]}=0 \qquad \text{ for $m>i-j$, $m\in\NN$,}
\end{equation}
and
\begin{equation}\label{C_prods}
c_{i,j}^{[i-j]}=\prod_{r=j}^{i -1}a_{r+1,r} \qquad \text{for $i>j$.}
\end{equation}
\end{Rem}
Next, we provide an inequality between $\eta_i$ and $\eta_j$, $i,j\in\{1,\dots,p\}$ (defined in \eqref{etadef}), see Barczy and Bezd\'any \cite[Lemma 4.2]{BarBez}.
\begin{Lem}\label{eta_ineq}
Let $p\in\NN$, and let $\bA\in\RR_+^{p\times p}$ be a lower triangular matrix such that $a_{i,i}=1$, $i\in\{1,\dots,p\}$.
If $c_{i,j}^{[m]}>0$ for some $i,j\in\{1,\dots,p\}$ and $m\in\{0,1,\ldots,p-1\}$,
 then $\eta_i\geq\eta_j+m$.
\end{Lem}

\section{Moments}\label{GWI_moments}

In this appendix, we recall some results from Appendix B in Barczy and Bezd\'any \cite{BarBez}, 
 and we also extend their Lemma B.4.
For the proofs of the results recalled from Appendix B in Barczy and Bezd\'any \cite{BarBez}, see their paper.
For the extended result, we give a detailed proof.

\begin{Lem}\label{Moments}
Let $(\bX_k)_{k\in\ZZ_+}$ be a $p$-type GWI process such that $\bX_0 = \bzero$ and the moment condition \eqref{baseconds} holds.
Then for each $k \in \NN$, we have $\EE(\bX_k \mid \cF_{k-1}^\bX) = \bA \bX_{k-1} + \bb$
 and
 \begin{gather}
  \EE(\bX_k)
   = \sum_{j=0}^{k-1} \bA^j \bb, \label{mean2} \\
  \var\bigl(\bX_k \mid \cF_{k-1}^\bX\bigr)
    = \var\bigl(\bM_k \mid \cF_{k-1}^\bX\bigr)
	= \EE(\bM_k \bM_k^\top \mid \cF_{k-1}^\bX)
     = \bV^{(0)} + \sum_{i=1}^p X_{k-1,i} \bV^{(i)}, \label{condVa2r}\\
   \var\bigl(\bX_k\bigr)
      = \sum_{j=0}^{k-1} \bA^j \EE(\bM_{k-j} \bM_{k-j}^\top) (\bA^\top)^j, \label{Var2}\\
  \EE(\bM_k \bM_k^\top)
  = \bV^{(0)} + \sum_{i=1}^p \EE(X_{k-1,i}) \bV^{(i)}. \label{Cov2}
 \end{gather}
\end{Lem}

Lemma \ref{Moments} is a special case of Lemma A.1 in Isp\'any and Pap \cite{IspPap2} for $p$-type GWI processes starting from $\bzero$.
For completeness, we note that Lemma A.1 in Isp\'any and Pap \cite{IspPap2} is stated only for critical
 $p$-type GWI processes, but its proof readily shows that it holds not only in the critical case.

The following lemma can be found in Barczy and Bezd\'any \cite[Lemma B.2]{BarBez}.

\begin{Lem}\label{EEX1}
Let $(\bX_k)_{k\in\ZZ_+}$ be a strongly critical $p$-type GWI process such that $\bX_0 = \bzero$,
 the moment condition \eqref{baseconds} holds, and suppose that the offspring mean matrix
 $\bA$ is lower triangular such that $a_{i,i}=1$, $i\in\{1,\dots,p\}$.
Then, for each $i,j\in\{1,\dots,p\}$, we have
 \begin{gather*}
  \EE(X_{k,i})=\OO(k^{\eta_i}), \qquad |\EE(M_{k,i}M_{k,j})|=\OO(k^{\eta_i\wedge\eta_j}),\qquad k\in\NN,
 \end{gather*}
 where $\eta_i$, $i\in\{1,\dots,p\}$, are defined in \eqref{etadef}.
\end{Lem}

The following lemma is an extension of Lemma B.4 in Barczy and Bezd\'any \cite{BarBez}.

\begin{Lem}\label{M_sum_growth}
Let $(\bX_k)_{k\in\ZZ_+}$ be a strongly critical $p$-type GWI process such that
 $\bX_0 = \bzero$, the moment condition \eqref{baseconds} holds, and suppose that the offspring mean matrix
 $\bA$ is lower triangular such that $a_{i,i}=1$, $i\in\{1,\dots,p\}$.
Then, for each $i\in\{1,\dots,p\}$, $m\in\ZZ_+$, and $T\in\RR_{++}$, we have
\begin{gather}\label{doobBound}
\EE\left(\sup_{t\in[0,T]}\left(\sum_{\ell=1}^{\nt}\binom{\nt-\ell}{m}(M_{\ell,i}+b_i)\right)^2\right)=\OO(n^{\eta_i+2m+1}),\qquad n\in\NN,
\end{gather}
where $\eta_i$, $i\in\{1,\dots,p\}$, are given in \eqref{etadef}.
\end{Lem}

\noindent\textbf{Proof.}
Let $T\in\RR_{++}$ and $i\in\{1,\ldots,p\}$ be fixed.
We proceed by induction on $m\in\ZZ_+$.
By the power mean inequality, we have
 \begin{equation*}
 \sup_{t\in[0,T]}\left(\sum_{\ell=1}^{\nt}(M_{\ell,i}+b_i)\right)^2\leq2\nT^2 b_i^2 + 2\sup_{t\in[0,T]}\left(\sum_{\ell=1}^\nt M_{\ell,i}\right)^2,
    \qquad n\in\NN.
 \end{equation*}
Using Doob's maximal inequality (see, e.g., Revuz and Yor \cite[Chapter II, Corollary (1.6)]{RevYor})
for the martingale $\left(\sum_{\ell=1}^k M_{\ell,i}\right)_{k\in\NN}$ (with respect to the filtration \ $(\cF_k^\bX)_{k\in\NN}$), we obtain
 \begin{equation*}
  \EE\left(\sup_{t\in[0,T]}\left(\sum_{\ell=1}^\nt M_{\ell,i}\right)^2\right)
     =\EE\left(\max_{k\in\{0,\dots,\nT\}}\left(\sum_{\ell=1}^k M_{\ell,i}\right)^2\right)
  	   \leq4\EE\left(\left(\sum_{\ell=1}^\nT M_{\ell,i}\right)^2\right),\qquad n\in\NN.
 \end{equation*}
Using that
 \begin{align*}
  \EE(M_{\ell,i}M_{j,i})
    & = \EE\big( \EE(M_{\ell,i}M_{j,i} \mid \cF_{j-1}^\bX) \big)
      = \EE\big( M_{\ell,i} \EE(M_{j,i} \mid \cF_{j-1}^\bX) \big)\\
    & = \EE\big( M_{\ell,i} \EE(X_{j,i}  - \EE(X_{j,i} \mid \cF_{j-1}^\bX ) \mid \cF_{j-1}^\bX) \big)
     = \EE(M_{\ell,i}\cdot 0)
     =0
 \end{align*}
 for each $n\in\NN$, $\ell\in\{ 1,\ldots,\nT-1\}$ and $j\in\{\ell+1,\ldots,\nT\}$, by Lemma \ref{EEX1}, we get that
 \begin{align*}
  \EE\left(\left(\sum_{\ell=1}^\nT M_{\ell,i}\right)^2\right)
	&=\EE\left( \sum_{\ell=1}^\nT M_{\ell,i}^2 + 2 \sum_{\ell=1}^{\nT-1} \sum_{j=\ell+1}^\nT M_{\ell,i} M_{j,i}\right)\\
	&=\sum_{\ell=1}^\nT\EE(M_{\ell,i}^2) + 2\sum_{\ell=1}^{\nT-1}\sum_{j=\ell+1}^\nT\EE(M_{\ell,i}M_{j,i})\\
	&=\sum_{\ell=1}^\nT\EE(M_{\ell,i}^2)=\sum_{\ell=1}^\nT\OO(\ell^{\eta_i})= \OO(n^{\eta_i+1}),
    \qquad n\in\NN.
\end{align*}
This together with $2\nT^2b_i^2=\OO(n^2)=\OO(n^{\eta_i+1})$, $n\in\NN$, since $\eta_i\geq1$, yields \eqref{doobBound} for the case $m=0$.

Now let $m\in\NN$ be fixed, and suppose that \eqref{doobBound} holds when $m$ is replaced by $0,\ldots,m-1$.
By Lemma \ref{specSums}, we have
\begin{equation}\label{help19}
\sum_{\ell=1}^{\nt}\binom{\nt-\ell}{m}(M_{\ell,i}+b_i)=n\int_{0}^{\frac{\nt}{n}}\sum_{\ell=1}^{\ns}\binom{\ns-\ell}{m-1}(M_{\ell,i}+b_i)\,\dd s, \qquad t\in\RR_+,\quad  n\in\NN.
\end{equation}
For each $n\in\NN$ and $t\in[0,n^{-1})$, we integrate on the set $\{0\}$ in the integral on the right hand side of \eqref{help19},
 so the integral is $0$ in this case.
Otherwise, for each $n\in\NN$ and $t\in[n^{-1},\infty)$, we have $\nt>0$, and, by (the measure-theoretic) Jensen's inequality,
 \begin{align*}
   \left(\int_{0}^{\frac{\nt}{n}}\sum_{\ell=1}^{\ns}\binom{\ns-\ell}{m-1}(M_{\ell,i}+b_i)\,\dd s\right)^2
   & =\frac{\nt^2}{n^2}\left(\frac{n}{\nt}\int_{0}^{\frac{\nt}{n}}\sum_{\ell=1}^{\ns}\binom{\ns-\ell}{m-1}(M_{\ell,i}+b_i)\,\dd s\right)^2\\
   & \leq \frac{\nt^2}{n^2}\cdot\frac{n}{\nt}\int_{0}^{\frac{\nt}{n}}\left(\sum_{\ell=1}^{\ns}\binom{\ns-\ell}{m-1}(M_{\ell,i}+b_i)\right)^2\,\dd s\\
   & \leq t\int_{0}^{t}\sup_{ r\in[0,t]}\left(\sum_{\ell=1}^{ \nr}\binom{\nr-\ell}{m-1}(M_{\ell,i}+b_i)\right)^2\,\dd s
 \end{align*}  
 \begin{align*}
   &\hspace{4.9cm} \leq t^2\sup_{ r\in[0,t]}\left(\sum_{\ell=1}^{ \nr}\binom{\nr-\ell}{m-1}(M_{\ell,i}+b_i)\right)^2.
 \end{align*}
Consequently, by \eqref{help19} and the induction hypothesis, we get that for all $T\in\RR_{++}$,
 \begin{align*}
 &\EE\left(\sup_{t\in[0,T]}\left(\sum_{\ell=1}^{\nt}\binom{\nt-\ell}{m}(M_{\ell,i}+b_i)\right)^2\right)\\
 &\qquad \leq n^2\EE\left(\sup_{t\in[0,T]}\left[ t^2\sup_{ r\in[0,t]}\left(\sum_{\ell=1}^{ \nr}\binom{\nr-\ell}{m-1}(M_{\ell,i}+b_i)\right)^2 \right] \right)\\
 &\qquad \leq n^2T^2\EE\left(\sup_{t\in[0,T]}\left(\sum_{\ell=1}^{\nt}\binom{\nt-\ell}{m-1}(M_{\ell,i}+b_i)\right)^2\right)\\
 &\qquad = n^2T^2 \OO(n^{\eta_i+2(m-1)+1}) = \OO(n^{\eta_i+2m+1}),\qquad n\in\NN,
 \end{align*}
 as desired.
\proofend

Finally, we formulate an easy consequence of Lemma \ref{M_sum_growth}.

\begin{Rem}\label{Rem_X2}
Let $(\bX_k)_{k\in\ZZ_+}$ be a $p$-type GWI process such that the conditions of
 Lemma \ref{M_sum_growth} are satisfied.
Recall that $q\in\{1,\dots,p\}$ is such that, for each $j\in\{1,\ldots,p\}$,
 we have that $\eta_j=1$ if and only if $j\leq q$,
with $\eta_i:=\eta_i(\bA)$, $i\in\{1,\dots,p\}$, defined in \eqref{etadef}.
Then  Lemma \ref{M_sum_growth} implies that $\EE(X_{k,i}^2)=\OO(k^2)$, $k\in\ZZ_+$, $i\in\{1,\dots,q\}$.
Indeed, for each $i\in\{1,\dots,q\}$, we have that $X_{k,i}=\sum_{\ell=1}^k(M_{\ell,i}+b_i)$, $k\in\ZZ_+$
 (following from Step 1/(c) in the proof of Theorem \ref{main_1} and \eqref{M_sums}),
 and hence
\[
\EE(X_{k,i}^2)
	=\EE(X_{\lfloor k\cdot 1\rfloor,i}^2)
	\leq\EE(\sup_{t\in[0,1]}X_{\lfloor kt\rfloor,i}^2)
	=\OO(k^{\eta_i+1})=\OO(k^2) , \qquad k\in\ZZ_+,
\]
 since $\eta_i=1$ for each $i\in\{1,\dots,q\}$.
\proofend
\end{Rem}

\section{A version of the continuous mapping theorem}
\label{CMT}

In this appendix, we recall several results from Appendix C in Barczy and Bezd\'any \cite{BarBez}, and we add new results as well.
Among others, we recall a version of the continuous mapping theorem
for $\RR^d$-valued stochastic processes with c\`adl\`ag paths,
 and then we present some auxiliary lemmas that make the application of
 this continuous mapping theorem easier throughout the proofs.
For the proofs of the results recalled from Appendix C in Barczy and Bezd\'any \cite{BarBez}, see their paper.
For the new results, we give detailed proofs.

A function $f : \RR_+ \to \RR^d$ is called c\`adl\`ag if it is right
 continuous with left limits.
Let $\DD(\RR_+, \RR^d)$ and $\CC(\RR_+, \RR^d)$ denote the space of
 all $\RR^d$-valued c\`adl\`ag and continuous functions on $\RR_+$, respectively.
Let $\cB(\DD(\RR_+, \RR^d))$ denote the Borel $\sigma$-algebra on
 $\DD(\RR_+, \RR^d)$ for the metric defined in Jacod and Shiryaev
 \cite[Chapter VI, (1.26)]{JacShi} (see also Billingsley \cite[Chapter 3, Section 16]{Bil}).
With this metric $\DD(\RR_+, \RR^d)$ is a
 complete and separable metric space and the topology induced by this metric is
 the so-called Skorokhod $J_1$ topology.
Recall that  $\skor{d}$ denotes the convergence in this Skorokhod
 topology for $\DD(\RR_+,\RR^d)$ (see Section \ref{Section_multi-type_GWI}).
Note that $\CC(\RR_+,\RR^d)\in\cB(\DD(\RR_+,\RR^d))$, see, e.g., Ethier and Kurtz \cite[Problem 3.11.25]{EthKur}.
For $\RR^d$-valued stochastic processes $(\bcY_t)_{t \in \RR_+}$ and
 $(\bcY^{(n)}_t)_{t \in \RR_+}$, $n \in \NN$, with c\`adl\`ag paths we write
 $\bcY^{(n)} \distr \bcY$ if the distribution of $\bcY^{(n)}$ on the
 space $(\DD(\RR_+,  \RR^d), \cB(\DD(\RR_+, \RR^d)))$ converges weakly to the
 distribution of $\bcY$ on the space
 $(\DD(\RR_+,  \RR^d), \cB(\DD(\RR_+, \RR^d)))$ as $n \to \infty$.
If $\xi$ and $\xi_n$, $n \in \NN$, are random elements with values in a metric space $(E,d)$,
 then we denote by $\xi_n \distr \xi$ the weak convergence of the
 distribution of $\xi_n$ on the space $(E, \cB(E))$ towards the
 distribution of $\xi$ on the space $(E, \cB(E))$ as $n \to \infty$,
 where $\cB(E)$ denotes the Borel $\sigma$-algebra on $E$ induced by
 the given metric $d$.

First, we establish some facts about measurability and continuity of mappings 
 (see Barczy and Bezd\'any \cite[Lemma C.4]{BarBez}).

\begin{Lem}\label{coordMeasurable}
Let $d,p,q\in\NN$, $\Phi_1:\DD(\RR_+,\RR^d)\to\DD(\RR_+,\RR^p)$ and $\Phi_2:\DD(\RR_+,\RR^d)\to\DD(\RR_+,\RR^q)$ be measurable mappings.
Then $\Phi:\DD(\RR_+,\RR^d)\to\DD(\RR_+,\RR^{p+q})$, $\Phi(f):=[\Phi_1(f),\Phi_2(f)]^\top$, $f\in\DD(\RR_+,\RR^d)$ is measurable as well.
\end{Lem}

The following lemma is a corollary of a result in Jacod and Shiryaev \cite[Chapter VI, Proposition 2.2 (b)]{JacShi}
 (see also Barczy and Bezd\'any \cite[Lemma C.5]{BarBez}).

\begin{Lem}\label{CoordConv}
Let $p,q\in\NN$, $f\in\DD(\RR_+,\RR^p)$, $f_n\in\DD(\RR_+,\RR^p)$, $n\in\NN$,
and $g\in\DD(\RR_+,\RR^q)$, $g_n\in\DD(\RR_+,\RR^q)$, $n\in\NN$, be
such that $f_n\skor{p} f$ as $n\to\infty$ and $g_n\skor{q} g$ as $n\to\infty$.
If $g\in\CC(\RR_+,\RR^q)$, $g_n\in\CC(\RR_+,\RR^q)$, $n\in\NN$,
then we have
\begin{equation*}
\begin{bmatrix}
f_n\\
g_n
\end{bmatrix}\skor{p+q}
\begin{bmatrix}
f\\
g
\end{bmatrix}\qquad \text{as $n\to\infty$.}
\end{equation*}
\end{Lem}

For the following corollary, see Barczy and Bezd\'any \cite[Lemma C.6]{BarBez}.

\begin{Cor}\label{ContMapVector}
Let $d,p,q\in\NN$, $\Phi_1:\DD(\RR_+,\RR^d)\to\DD(\RR_+,\RR^p)$ and $\Phi_2:\DD(\RR_+,\RR^d)\to\DD(\RR_+,\RR^q)$ be continuous mappings.
If $\Phi_1(\DD(\RR_+,\RR^d))\subset\CC(\RR_+,\RR^p)$ and $\Phi_2(\DD(\RR_+,\RR^d))\subset\CC(\RR_+,\RR^q)$,
then $\Phi:\DD(\RR_+,\RR^d)\to\DD(\RR_+,\RR^{p+q})$, $\Phi(f):=[\Phi_1(f),\Phi_2(f)]^\top$, $f\in\DD(\RR_+,\RR^d)$ is continuous as well, and $\Phi(\DD(\RR_+,\RR^d))\subset\CC(\RR_+,\RR^{p+q})$.
\end{Cor}

The following version of the continuous mapping theorem can be found for
 example in Theorem 3.27 of Kallenberg \cite{Kal}.

\begin{Lem}\label{lemma:kallenberg}
Let $(S, d_S)$ and $(T, d_T)$ be metric spaces and
 $(\xi_n)_{n \in \NN}$ and $\xi$ be random elements with values in $S$
 such that $\xi_n \distr \xi$ as $n \to \infty$.
Let $f : S \to T$ and $f_n : S \to T$, $n \in \NN$, be measurable
 mappings and $C \in \cB(S)$ such that $\PP(\xi \in C) = 1$ and
 $\lim_{n \to \infty} d_T(f_n(s_n), f(s)) = 0$ if
 $\lim_{n \to \infty} d_S(s_n,s) = 0$ and $s \in C$, $s_n\in S$, $n\in\NN$.
Then $f_n(\xi_n) \distr f(\xi)$ as $n \to \infty$.
\end{Lem}

For functions $f$ and $f_n$, $n \in \NN$, in $\DD(\RR_+, \RR^d)$,
 we write $f_n \lu f$ if $(f_n)_{n \in \NN}$ converges to $f$
 locally uniformly, i.e., if $\sup_{t \in [0,T]} \|f_n(t) - f(t)\| \to 0$
 as $n \to \infty$ for all $T \in \RR_{++}$.
For measurable mappings $\Phi : \DD(\RR_+, \RR^d) \to \DD(\RR_+, \RR^q)$
 and $\Phi_n : \DD(\RR_+, \RR^d) \to \DD(\RR_+, \RR^q)$, $n \in \NN$,
 we will denote by $C_{\Phi, (\Phi_n)_{n \in \NN}}$ the set of all functions
 $f \in \CC(\RR_+, \RR^d)$ \ for which \ $\Phi_n(f_n) \lu \Phi(f)$ whenever
 $f_n \lu f$ \ with \ $f_n \in \DD(\RR_+, \RR^d)$, $n \in \NN$.

One can formulate the following consequence of Lemma \ref{lemma:kallenberg}.

\begin{Lem}\label{Conv2Funct}
Let $d,\,q\in\NN$.
Let $(\bcU_t)_{t \in \RR_+}$ and $(\bcU^{(n)}_t)_{t \in \RR_+}$, $n \in \NN$,
 be $\RR^d$-valued stochastic processes with c\`adl\`ag paths such that
 $\bcU^{(n)} \distr \bcU$ as $n\to\infty$.
Let $\Phi : \DD(\RR_+, \RR^d) \to \DD(\RR_+, \RR^q)$ and
  $\Phi_n : \DD(\RR_+, \RR^d) \to \DD(\RR_+, \RR^q)$, $n \in \NN$, be
 measurable mappings such that there exists
 $C \subset C_{\Phi,(\Phi_n)_{n\in\NN}}$ with $C \in \cB(\DD(\RR_+, \RR^d))$ and $\PP(\bcU \in C) = 1$.
Then $\Phi_n(\bcU^{(n)}) \distr \Phi(\bcU)$ as $n\to\infty$.
\end{Lem}

\begin{Lem}\label{ConvIntersect}
Let $d,q,p\in\NN$, and let $\Phi^{(1)}$ and $\Phi_n^{(1)}$, $n\in\NN$,
 be measurable mappings from $\DD(\RR_+,\RR^d)$ to $\DD(\RR_+,\RR^p)$,
 and $\Phi^{(2)}$ and $\Phi_n^{(2)}$, $n\in\NN$, be measurable mappings from $\DD(\RR_+,\RR^d)$ to $\DD(\RR_+,\RR^q)$.
Then $\Phi:=[\Phi^{(1)}, \Phi^{(2)}]^\top$, and $\Phi_n:=[\Phi_n^{(1)}, \Phi_n^{(2)}]^\top$, $n\in\NN$ are measurable mappings from $\DD(\RR_+,\RR^d)$ to $\DD(\RR_+,\RR^{p+q})$,
and we have that $ C_{\Phi^{(1)},(\Phi_n^{(1)})_{n\in\NN}}\cap C_{\Phi^{(2)},(\Phi_n^{(2)})_{n\in\NN}}\subset C_{\Phi,(\Phi_n)_{n\in\NN}}$.
\end{Lem}

\noindent\textbf{Proof.}
The measurability of the given maps follows from Lemma \ref{coordMeasurable}.
Let $f_n \in\DD(\RR_+,\RR^d)$, $n\in\NN$, and $f\in C_{\Phi^{(1)},(\Phi_n^{(1)})_{n\in\NN}}\cap C_{\Phi^{(2)},(\Phi_n^{(2)})_{n\in\NN}}$
 be such that $f_n\lu f$ as $n\to\infty$.
Then for each $n\in\NN$ and $T\in\RR_{++}$, we have
\begin{align*}
&\sup_{t\in[0,T]}\| (\Phi_n(f_n))(t)  -  (\Phi(f))(t)\|
 =\sup_{t\in[0,T]}\left\|
  \begin{bmatrix}
   (\Phi_n^{(1)}(f_n))(t)\\
   (\Phi_n^{(2)}(f_n))(t)
  \end{bmatrix}
 -\begin{bmatrix}
   (\Phi^{(1)}(f))(t)\\
   (\Phi^{(2)}(f))(t)
  \end{bmatrix}\right\|\\
&\qquad =\sup_{t\in[0,T]}\left\|
\begin{bmatrix}
(\Phi_n^{(1)}(f_n))(t) - (\Phi^{(1)}(f))(t)\\
0
\end{bmatrix}
+\begin{bmatrix}
0\\
(\Phi_n^{(2)}(f_n))(t) - (\Phi^{(2)}(f))(t)
\end{bmatrix}\right\|\\
&\qquad \leq\sup_{t\in[0,T]}\left\|
\begin{bmatrix}
(\Phi_n^{(1)}(f_n))(t) - (\Phi^{(1)}(f))(t)\\
0
\end{bmatrix}\right\|
+\sup_{t\in[0,T]}\left\|\begin{bmatrix}
0\\
(\Phi_n^{(2)}(f_n))(t) - (\Phi^{(2)}(f))(t)
\end{bmatrix}\right\|\\
&\qquad=\sup_{t\in[0,T]}\left\|
(\Phi_n^{(1)}(f_n))(t) - (\Phi^{(1)}(f))(t)
\right\|
+\sup_{t\in[0,T]}\left\|
(\Phi_n^{(2)}(f_n))(t) - (\Phi^{(2)}(f))(t)
\right\|.
\end{align*}
Since $f_n\lu f$ as $n\to\infty$ and $f\in C_{\Phi^{(1)},(\Phi_n^{(1)})_{n\in\NN}}\cap C_{\Phi^{(2)},(\Phi_n^{(2)})_{n\in\NN}}$,
we have that both terms on the right converge to $0$ as $n\to\infty$.
Thus $\Phi_n(f_n)\lu\Phi(f)$ as $n\to\infty$,
meaning that $f\in C_{\Phi,(\Phi_n)_{n\in\NN}}$, which implies the desired assertion.
\proofend

Let us introduce some mappings on the space $\DD(\RR_+,\RR^d)$, where $d\in\NN$.

\begin{Def}\label{maps}
Let $d,q\in\NN$. For $\bB\in\RR^{q\times d}$, $m\in\ZZ_+$, $n\in\NN$, and $s\in\RR_+$, we define the following maps:
\begin{itemize}
\item[(i)] $\mathrm{Id}_{d}:\DD(\RR_+,\RR^d)\to\DD(\RR_+,\RR^d)$, $(\mathrm{Id}_{d}(f))(t):=f(t)$ for $f\in\DD(\RR_+,\RR^d)$ and $t\in\RR_+$.
\item[(ii)] $\Psi_\bB:\DD(\RR_+,\RR^d)\to\DD(\RR_+,\RR^q)$, $(\Psi_\bB(f))(t):=\bB f(t)$ for $f\in\DD(\RR_+,\RR^d)$ and $t\in\RR_+$.
\item[(iii)] $\pi_s:\DD(\RR_+,\RR^d)\to\RR^d$, $\pi_s(f):=f(s)$ for $f\in\DD(\RR_+,\RR^d)$.
\item[(iv)] $\tau_n:\DD(\RR_+,\RR^d)\to\DD(\RR_+,\RR^d)$, $(\tau_n(f))(t):=f(\nt/n)$ for $f\in\DD(\RR_+,\RR^d)$ and $t\in\RR_+$.
\item[(v)] $\varphi:\DD(\RR_+,\RR^d)\to\DD(\RR_+,\RR^d)$, $(\varphi(f))(t):=\int_0^tf(s)\,\dd s$ for $f\in\DD(\RR_+,\RR^d)$ and $t\in\RR_+$.
\item[(vi)] $\varphi_n:\DD(\RR_+,\RR^d)\to\DD(\RR_+,\RR^d)$, $\varphi_n:=\tau_n\circ\varphi$.
\item[(vii)] $\Phi^{(m)}:\DD(\RR_+,\RR^d)\to\DD(\RR_+,\RR^{(m+1)d})$ and $\Phi_n^{(m)}:\DD(\RR_+,\RR^d)\to\DD(\RR_+,\RR^{(m+1)d})$,
\begin{gather}\label{Phi_def}
(\Phi^{(m)}(f))(t):=\begin{bmatrix}
                      (\varphi^{\circ (0)}(f))(t) \\
                    (\varphi^{ \circ (1)}(f))(t)\\
		\vdots\\
		(\varphi^{\circ (m)}(f))(t)
                  \end{bmatrix},	\qquad
(\Phi_n^{(m)}(f))(t):=\begin{bmatrix}
                      (\varphi_n^{\circ (0)}(f))(t) \\
                    (\varphi_n^{ \circ (1)}(f))(t)\\
		\vdots\\
		(\varphi_n^{\circ (m)}(f))(t)
                  \end{bmatrix}
\end{gather}
for $f\in\DD(\RR_+,\RR^d)$ and $t\in\RR_+$, where we recall the notation for $m$-fold compositions:
 $\varphi^{\circ(0)}$ and $\varphi_n^{\circ(0)}$ are the identity map $\mathrm{Id}_d$ on $\DD(\RR_+,\RR^d)$,
 and $\varphi^{\circ(m+1)}:=\varphi\circ\varphi^{\circ(m)}$, $\varphi_n^{\circ(m+1)}:=\varphi_n\circ\varphi_n^{\circ(m)}$ for $m\in\ZZ_+$.
\end{itemize}
\end{Def}

 Next, we present some auxiliary lemmas which are used throughout the proofs.
The following Lemma \ref{ConvMatrixMult} can be found, e.g., 
 in Barczy and Bezd\'any \cite[Lemma C.3]{BarBez}.

\begin{Lem}\label{ConvMatrixMult}
Let $d,q\in\NN$, $\bB\in\RR^{q\times d}$, and let $\Psi_\bB$ be the map given in part (ii) of Definition \ref{maps}.
Then $ \Psi_\bB$ is continuous  (in particular, measurable).
Further, if $(\bcU_t)_{t \in \RR_+}$ and $(\bcU^{(n)}_t)_{t \in \RR_+}$, $n \in \NN$,
 are $\RR^d$-valued stochastic processes with c\`adl\`ag paths such that
 $\bcU^{(n)} \distr \bcU$ as $n\to\infty$, then $\bB\bcU^{(n)} \distr \bB\bcU$ as $n\to\infty$.
\end{Lem}

\begin{Lem}\label{time_transform}
Let $d\in\NN$, and for each $n\in\NN$, let $\tau_n$ be the map given in part (iv) of Definition \ref{maps}.
We have that for each $n\in\NN$, $\tau_n$ is Borel measurable.
Further, if $(\bcU_t)_{t \in \RR_+}$ and $(\bcU^{(n)}_t)_{t \in \RR_+}$, $n \in \NN$,
 are $\RR^d$-valued stochastic processes with c\`adl\`ag paths such that $\PP(\bcU\in\CC(\RR_+,\RR^d))=1$ and
 $\bcU^{(n)} \distr \bcU$ as $n\to\infty$, then we have that $\tau_n(\bcU^{(n)}) \distr \bcU$ as $n\to\infty$.
\end{Lem}

\noindent\textbf{Proof.}
Let $n\in\NN$ be fixed for the rest of this paragraph.
Using that the mappings $\pi_t$, $t\in\RR_+$, generate the Borel $\sigma$-algebra on $\DD(\RR_+, \RR^d)$
(see, e.g., Jacod and Shiryaev \cite[Chapter VI, Theorem 1.14, part c)]{JacShi} or Billingsley \cite[part (iii) of Theorem 16.6]{Bil}),
 to check the Borel measurability of $\tau_n$, it is enough to verify
 that the mappings $\pi_t\circ \tau_n:\DD(\RR_+,\RR^d)\to\RR^d$
 are measurable for  all $t\in\RR_+$.
Indeed, it is sufficient to check that $\tau_n^{-1}(\pi_t^{-1}(B))\in\cB(\DD(\RR_+,\RR^d))$ for all $t\in\RR_+$
 and $B\in\cB(\RR^d)$ (see, e.g., Klenke \cite[Theorem 1.81]{Kle}) and we have that
 $\tau_n^{-1}(\pi_t^{-1}(B)) = (\pi_t\circ \tau_n)^{-1}(B)$, $t\in\RR_+$, $B\in\cB(\RR^d)$.
Notice that $\pi_t\circ \tau_n=\pi_{\nt/n}$ for all $t\in\RR_+$.
Since the natural projection $\pi_{\nt/n}$ onto $\frac{\nt}{n}$ is measurable for all $t\in\RR_+$
 (see, e.g., Billingsley \cite[Theorem 16.6, part (i)]{Bil}),
we get that $\pi_t \circ \tau_n$ is measurable as well, and thus so is $\tau_n$.

Let $f_n\in\DD(\RR_+,\RR^d)$, $n\in\NN$, and $f\in\CC(\RR_+,\RR^d)$ be such that $f_n\lu f$ as $n\to\infty$.
Then for each $n\in\NN$ and $T\in\RR_{++}$, we have
\begin{align*}
\sup_{t\in[0,T]}\|(\tau_n(f_n))(t)-f(t)\|
	\leq&\sup_{t\in[0,T]}\left\|f_n\left(\frac{\nt}{n}\right)-f\left(\frac{\nt}{n}\right)\right\|+\sup_{t\in[0,T]}\left\|f\left(\frac{\nt}{n}\right)-f(t)\right\|\\
	=&\sup_{t\in[0,T]}\left\|(\tau_n(f_n-f))(t)\right\|+\sup_{t\in[0,T]}\left\|f\left(\frac{\nt}{n}\right)-f(t)\right\|\\
	=&\max_{k\in\{0,\dots,\nT\}}\left\|f_n\left(\frac{k}{n}\right)-f\left(\frac{k}{n}\right)\right\|+\sup_{t\in[0,T]}\left\|f\left(\frac{\nt}{n}\right)-f(t)\right\|\\
	\leq&\sup_{t\in[0,T]}\|f_n(t)-f(t)\|+\max_{k\in\{0,\dots,\nT\}}\sup_{t\in\left[\frac{k}{n},\frac{k+1}{n}\right)}\left\|f\left(\frac{k}{n}\right)-f(t)\right\|.
\end{align*}
Since $f_n\lu f$ as $n\to\infty$, the first term tends to $0$ as $n\to\infty$.
Further, using that $[0,T]$, $T\in\RR_{++}$, is compact and $f$ is continuous, we have that
 $f$ is uniformly continuous on $[0,T]$, and thus the second term also tends to $0$,
and we get that $\tau_n(f_n)\lu f$ as $n\to\infty$.
Consequently, we have $C_{\mathrm{Id}_d,(\tau_n)_{n\in\NN}}=\CC(\RR_+,\RR^d)$,
thus by Lemma \ref{Conv2Funct} (with the choice of $C=\CC(\RR_+,\RR^d))$
we have that  $\tau_n(\bcU^{(n)}) \distr \bcU$ as $n\to\infty$.
\proofend

The following Lemma \ref{IntegralCont} can be found, e.g., 
 in Barczy and Bezd\'any \cite[Lemma C.7]{BarBez}.

\begin{Lem}\label{IntegralCont}
For each $d\in\NN$, the map $\varphi$ given in part (v) of Definition \ref{maps}
 is continuous, and maps $\DD(\RR_+,\RR^d)$ into $\CC(\RR_+,\RR^d)$.
\end{Lem}

\begin{Lem}\label{ConvIteratedIntegral}
Let $d\in\NN$ and $m\in\ZZ_+$. The maps $\Phi^{(m)}$ and $\Phi_n^{(m)}$, $n\in\NN$, defined in part (vii) of Definition \ref{maps}
are measurable, and if $(\bcU_t)_{t \in \RR_+}$ and $(\bcU^{(n)}_t)_{t \in \RR_+}$, $n \in \NN$,
 are $\RR^d$-valued stochastic processes with c\`adl\`ag paths such that $\PP(\bcU\in\CC(\RR_+,\RR^d))=1$ and
 $\bcU^{(n)} \distr \bcU$ as $n\to\infty$, then we have that $\Phi_n^{(m)}(\bcU^{(n)}) \distr \Phi^{(m)}(\bcU)$ as $n\to\infty$.
\end{Lem}

\noindent\textbf{Proof.}
Let $n\in\NN$ be fixed.
By Lemma \ref{coordMeasurable}, to check the measurability of $\Phi^{(m)}$ and $\Phi_n^{(m)}$, it is enough to check the measurability
of their coordinates in \eqref{Phi_def}. The measurability of the identity map is clear, and $\varphi$ is measurable by Lemma \ref{IntegralCont}.
Since $\varphi_n=\tau_n\circ\varphi$, and $\tau_n$ is (Borel) measurable by Lemma \ref{time_transform},
we have the measurability of $\varphi_n$ as well, as Borel measurability is preserved by composition.
Using again that taking composition preserves Borel measurability, we have that the maps $\varphi^{\circ(m)}$ and $\varphi_n^{\circ(m)}$
are Borel measurable for each $m\in\ZZ_+$, implying the measurability of $\Phi^{(m)}$ and $\Phi_n^{(m)}$ due to Lemma \ref{coordMeasurable}.

Next, we show that $C_{\Phi^{(m)},(\Phi_n^{(m)})_{n\in\NN}}=\CC(\RR_+,\RR^d)$ for each $m\in\ZZ_+$.
First, note that  $C_{\varphi^{\circ (0)}, (\varphi_n^{ \circ (0)})_{n\in\NN}}=\CC(\RR_+,\RR^d)$ (due to the continuity of the identity map $\mathrm{Id}_d$).
Let $f\in\CC(\RR_+,\RR^d)$ and $f_n\in \DD(\RR_+,\RR^d)$, $n\in\NN$, be such that $f_n\lu f$ as $n\to\infty$.
For all $n\in\NN$ and $t\in\RR_+$, we have
\begin{align*}
\Vert (\varphi_n(f_n))(t) - (\varphi(f))(t)\Vert
	=\left\Vert \int_0^{\frac{\nt}{n}}f_n(s)\,\dd s-\int_0^tf(s)\,\dd s\right\Vert,
\end{align*}
and thus
\begin{align*}
\Vert (\varphi_n(f_n))(t) - (\varphi(f))(t)\Vert
	\leq\int_0^{\frac{\nt}{n}}\left\Vert f_n(s)-f(s)\right\Vert \,\dd s+\int_{\frac{\nt}{n}}^t\left\Vert f(s)\right\Vert \,\dd s.
\end{align*}
Using also that
 \[
   \sup_{t\in[0,T]} \left( \frac{nt-\nt}{n} \right)
      \leq \frac{1}{n}, \qquad T\in\RR_{++}, \; n\in\NN,
 \]
 for all $T\in\RR_{++}$ and $n\in\NN$, we get that
\begin{align*}
\sup_{t\in[0,T]}\Vert (\varphi_n(f_n))(t) - (\varphi(f))(t)\Vert
	\leq& \frac{\nT}{n}\sup_{t\in[0,T]}\left\Vert f_n(t)-f(t)\right\Vert
            +  \sup_{t\in[0,T]} \left( \frac{nt-\nt}{n} \right)  \sup_{t\in[0,T]}\left\Vert f(t)\right\Vert\\
	\leq& T\sup_{t\in[0,T]}\left\Vert f_n(t)-f(t)\right\Vert +n^{-1}\sup_{t\in[0,T]}\left\Vert f(t)\right\Vert .
\end{align*}
For all $T\in\RR_{++}$, we have $\sup_{t\in[0,T]}\left\Vert f(t)\right\Vert<\infty$, since $f\in\CC(\RR_+,\RR^d)$,
thus $n^{-1}\sup_{t\in[0,T]}\left\Vert f(t)\right\Vert\to0$ as $n\to\infty$.
Since $f_n\lu f$ as $n\to\infty$, we have $\sup_{t\in[0,T]}\left\Vert f_n(t)-f(t)\right\Vert\to0$ as $n\to\infty$ for all $T\in\RR_+$.
Thus we get that $\varphi_n(f_n)\lu\varphi(f)$ as $n\to\infty$,
and so $C_{\varphi,(\varphi_n)_{n\in\NN}}=\CC(\RR_+,\RR^d)$.

Now let $m\in\NN\setminus\{1\}$ be fixed, and suppose that
$C_{\varphi^{\circ (j)},(\varphi_n^{\circ (j)})_{n\in\NN}}=\CC(\RR_+,\RR^d)$ holds for $j\in\{1,\ldots,m-1\}$.
Since $C_{\varphi^{\circ (m-1)},(\varphi_n^{\circ (m-1)})_{n\in\NN}}=\CC(\RR_+,\RR^d)$,
we have that $\varphi_n^{ \circ (m-1)}(f_n)\lu\varphi^{ \circ (m-1)}(f)$ as $n\to\infty$.
Then, since $\varphi^{\circ (m-1)}(f)=\varphi(\varphi^{\circ (m-2)}(f))\in \CC(\RR_+,\RR^d)$ (due to Lemma \ref{IntegralCont})
  and $C_{\varphi,(\varphi_n)_{n\in\NN}}=\CC(\RR_+,\RR^d)$, we have that $\varphi^{\circ (m-1)}(f)\in C_{\varphi,(\varphi_n)_{n\in\NN}}$.
Using the definition of $C_{\varphi,(\varphi_n)_{n\in\NN}}$ and that $\varphi_n^{\circ  (m-1)}(f_n)\in\DD(\RR_+,\RR^d)$, $n\in\NN$, we have that
 $\varphi_n(\varphi_n^{\circ (m-1)}(f_n))\lu\varphi(\varphi^{\circ (m-1)}(f))$ as $n\to\infty$.
By the definition of composition of functions, this means that $\varphi_n^{\circ (m)}(f_n)\lu\varphi^{\circ (m)}(f)$ as $n\to\infty$,
 thus $C_{\varphi^{\circ (m)},(\varphi_n^{\circ (m)})_{n\in\NN}}=\CC(\RR_+,\RR^d)$.

Using  Lemma \ref{ConvIntersect} and the fact that $C_{\varphi^{\circ (0)}, (\varphi_n^{ \circ (0)})_{n\in\NN}}=\CC(\RR_+,\RR^d)$, we have that
 \[
 \bigcap_{r=0}^{m}C_{\varphi^{\circ (r)},(\varphi_n^{\circ (r)})_{ n\in\NN }}=\bigcap_{r=0}^{m}\CC(\RR_+,\RR^d)
     = \CC(\RR_+,\RR^d) \subset C_{\Phi^{(m)},(\Phi_n^{(m)})_{n\in\NN}},
 \]
 and since $C_{\Phi^{(m)},(\Phi_n^{(m)})_{n\in\NN}}\subset\CC(\RR_+,\RR^d)$,
 we have $C_{\Phi^{(m)},(\Phi_n^{(m)})_{n\in\NN}}=\CC(\RR_+,\RR^d)$.

Consequently, one can apply Lemma \ref{Conv2Funct} by choosing $C:=\CC(\RR_+,\RR^d)$,
 and then, for each $m\in\ZZ_+$, we have that $\Phi_n^{(m)}(\bcU^{(n)}) \distr \Phi^{(m)}(\bcU)$ as $n\to\infty$, as desired.
\proofend

Finally, we recall a result on weak convergence for the sum of stochastic processes with c\`adl\`ag paths 
 due to Jacod and Shiryaev \cite[Lemma VI.3.31]{JacShi}.
  
\begin{Lem}\label{Lem_JacShi}
Let $(\bcY^{(n)}_t)_{t\in\RR_+}$, $n\in\NN$,  $(\bcY_t)_{t\in\RR_+}$, and  $(\bcZ^{(n)}_t)_{t\in\RR_+}$, $n\in\NN$,
 be $\RR^d$-valued stochastic processes with c\`adl\`ag paths on a probability space $(\Omega,\cA,\PP)$.
Suppose that $\bcY^{(n)}\distr \bcY$ as $n\to\infty$ and 
 \[
   \lim_{n\to\infty} \PP\Big( \sup_{t\in[0,T]} \Vert \bcZ^{(n)}_t\Vert > \vare \Big) = 0
     \qquad \text{for all \ $T\in\RR_{++}$ \ and \ $\vare\in\RR_{++}$.}
 \] 
Then $\bcY^{(n)} + \bcZ^{(n)} \distr \bcY$ as $n\to\infty$. 
\end{Lem}

\section{Convergence of random step processes}
\label{section_conv_step_processes}

We recall a result about convergence of random step processes towards a diffusion process, see Isp\'any and Pap \cite{IspPap}.

\begin{Thm}\label{Conv2DiffThm}
Let $\bbeta : \RR_+ \times \RR^d \to \RR^d$ and $\bgamma : \RR_+ \times \RR^d \to \RR^{d \times r}$
 be continuous functions.
Assume that uniqueness in the sense of probability law holds for the SDE
 \begin{equation}\label{SDE}
  \dd \, \bcU_t
  = \bbeta (t, \bcU_t) \, \dd t + \bgamma (t, \bcU_t) \, \dd \bcW_t ,
  \qquad t \in \RR_+,
 \end{equation}
 with initial value $\bcU_0 = \bu_0$ for all $\bu_0 \in \RR^d$, where
 $(\bcW_t)_{t\in\RR_+}$ is an $r$-dimensional standard Wiener process.
Let $(\bcU_t)_{t\in\RR_+}$ be a solution of \eqref{SDE} with initial value
 $\bcU_0 = \bzero \in \RR^d$.

For each $n \in \NN$, let $(\bU^{(n)}_k)_{k\in\ZZ_+}$ be a sequence of $d$-dimensional random vectors adapted to a filtration
  $(\cF^{(n)}_k)_{k\in\ZZ_+}$ (i.e., $\bU^{(n)}_k$ is $\cF^{(n)}_k$-measurable) such that
  $\EE(\Vert \bU^{(n)}_k\Vert^2)<\infty$ for each $n,k \in \NN$.
Let
 \[
   \bcU^{(n)}_t := \sum_{k=0}^{\nt}  \bU^{(n)}_k \, ,
   \qquad t \in \RR_+, \quad n \in \NN .
 \]
Suppose that $\bcU^{(n)}_0 = \bU^{(n)}_0 \distr \bzero$ as $n\to\infty$ and that for all $T \in \RR_{++}$,
 \begin{enumerate}
  \item[\textup{(i)}]
   $\sup\limits_{t\in[0,T]}
     \biggl\|\sum\limits_{k=1}^{\nt}
              \EE\bigl(\bU^{(n)}_k \mid \cF^{(n)}_{k-1}\bigr)
             - \int_0^t \bbeta(s,\bcU^{(n)}_s) \dd s\biggr\|
   \stoch 0$ \ as \ $n\to\infty$,  \\
  \item[\textup{(ii)}]
   $\sup\limits_{t\in[0,T]}
     \biggl\|\sum\limits_{k=1}^{\nt}
              \var\bigl(\bU^{(n)}_k \mid \cF^{(n)}_{k-1}\bigr)
             - \int_0^t
                \bgamma(s,\bcU^{(n)}_s) \bgamma(s,\bcU^{(n)}_s)^\top
                \dd s\biggr\|
   \stoch 0$ \ as \ $n\to\infty$, \\
  \item[\textup{(iii)}]
   $\sum\limits_{k=1}^{\lfloor nT \rfloor}
     \EE\bigl(\|\bU^{(n)}_k\|^2 \bbone_{\{\|\bU^{(n)}_k\| > \theta\}}
             \bmid \cF^{(n)}_{k-1}\bigr)
   \stoch 0$ \ as \ $n\to\infty$ \ for all \ $\theta \in \RR_{++}$.
 \end{enumerate}
Then $\bcU^{(n)} \distr \bcU$ as $n \to \infty$.
\end{Thm}

Note that in (ii) of Theorem \ref{Conv2DiffThm}, $\|\cdot\|$ denotes
 an operator norm, while in (i) it denotes a vector norm.

The following result is about the asymptotic behavior of a single-type Galton-Watson process with immigration in the critical case,
 it is due to Wei and Winnicki \cite[Theorem 2.1]{WW}.

\begin{Thm}\label{thm:critical}
Let $(X_k)_{k\in\ZZ_+}$ be a single-type Galton-Watson process with immigration such that
  $\EE(\xi^2) < \infty$, $\EE(\vare^2) < \infty$, $\EE(\xi) = 1$ (critical case) and $\EE(X_0^2) < \infty$.
Then
 \begin{equation*}
  (n^{-1} X_\nt)_{t\in\RR_+} \distr (\cX_t)_{t\in\RR_+} \qquad
  \text{as \ $n \to \infty$,}
 \end{equation*}
 where the limit process $(\cX_t)_{t\in\RR_+}$ is the pathwise unique strong solution of the SDE
 \begin{equation*}
  \dd \cX_t
  = \EE(\vare) \, \dd t
    + \sqrt{\var(\xi) \, \cX_t^+} \, \dd \cW_t , \qquad t\in\RR_+,
  \qquad \cX_0 = 0 ,
 \end{equation*}
 where $(\cW_t)_{t\in\RR_+}$ is a standard Wiener process.
\end{Thm}


\section*{Data availability statements}
Data sharing is not applicable to this article as no datasets were generated or analyzed during the current study.

\section*{Declaration of competing interest}
The authors declare that they have no conflicts of interest.

\end{document}